%% file: vogt.tex
\documentclass[graybox]{svmult}

\usepackage{type1cm}        
\usepackage{makeidx}         
\usepackage{graphicx}        
\usepackage{multicol}        
\usepackage[bottom]{footmisc}

\usepackage{newtxtext}       %
\usepackage{newtxmath}       

\makeindex             

\input{header.tex}

\begin{document}

\title*{Lifting methods for manifold-valued variational problems}
\author{Thomas Vogt, Evgeny Strekalovskiy, Daniel Cremers, Jan Lellmann}

\institute{%
Thomas Vogt, Jan Lellmann \at
Institute of Mathematics and Image Computing,
University of L\"ubeck,
Maria-Goeppert-Str. 3,
23562 L\"ubeck,
Germany,
\email{vogt@mic.uni-luebeck.de}, \email{lellmann@mic.uni-luebeck.de} %
\and %
Evgeny Strekalovskiy\at
Technical University Munich,
85748 Garching,
Germany.\at
Now at Google Germany GmbH, 
\email{evgeny.strekalovskiy@gmail.com}
\and %
Daniel Cremers\at
Technical University Munich,
85748 Garching,
Germany,
\email{cremers@tum.de}
}

\maketitle

\abstract{Lifting methods allow to transform hard variational problems such as segmentation and optical flow estimation into convex problems in a suitable higher-dimensional space. The lifted models can then be efficiently solved to a global optimum, which allows to find approximate global minimizers of the original problem. Recently, these techniques have also been applied to problems with values in a manifold. We provide a review of such methods in a refined framework based on a finite element discretization of the range, which extends the concept of sublabel-accurate lifting to manifolds. We also generalize existing methods for total variation regularization to support general convex regularization. 
}

\input{chpt/intro.tex}

\input{chpt/sec1.tex}

\input{chpt/sec2.tex}

\input{chpt/sec3.tex}
%

\begin{acknowledgement}
The authors acknowledge support through DFG grant LE \mbox{4064/1-1}
``Functional Lifting 2.0: Efficient Convexifications for Imaging and Vision''
and NVIDIA Corporation.
\end{acknowledgement}

\bibliographystyle{spmpsci}
\bibliography{vogt}

\end{document}

%% file: header.tex
\usepackage{tikz}
\usetikzlibrary{intersections}
\usetikzlibrary{arrows}

\usetikzlibrary{positioning}
\usetikzlibrary{calc}
\usetikzlibrary{spy}
\newcommand*{\ExtractCoordinate}[1]{\path (#1); \pgfgetlastxy{\XCoord}{\YCoord}; }%
\newcommand*{\ExtractImgDims}[1]{
    \ExtractCoordinate{$(#1.south west)$};
    \pgfmathsetmacro{\imgx}{\XCoord}
    \pgfmathsetmacro{\imgy}{\YCoord}
    \ExtractCoordinate{$(#1.north east)$};
    \pgfmathsetmacro{\imgw}{\XCoord - \imgx}
    \pgfmathsetmacro{\imgh}{\YCoord - \imgy}
}
\newcommand*{\RelativeSpy}[4]{
    \ExtractImgDims{#2};
    \begin{scope}[x=\imgw,y=\imgh,xshift=\imgx,yshift=\imgy]
        \coordinate (spyroi-#1) at #3;
        \coordinate (spypos-#1) at #4;
        \spy[anchor=center] on (spyroi-#1) in node[anchor=center] at (spypos-#1);
    \end{scope}
}

\DeclareMathOperator{\TV}{TV}
\DeclareMathOperator{\kernel}{ker} 
\DeclareMathOperator{\Div}{div}
\DeclareMathOperator*{\argmin}{arg\,min}
\newcommand{\R}{\mathbb{R}}
\newcommand{\IP}{\mathbb{P}}
\newcommand{\IS}{\mathbb{S}}

\newcommand{\IK}{\mathcal{K}}
\newcommand{\IM}{\mathcal{M}}
\newcommand{\ITh}{\mathcal{T}_h}
\newcommand{\IMh}{\IM_h}
\newcommand{\IfM}{\mathfrak{M}}

\usepackage{soulutf8}
\soulregister\cite7
\soulregister\ref7
\soulregister\pageref7
\soulregister\eqref7

%% file: chpt/intro.tex
\section{Introduction}

Consider a variational image processing or general data analysis problem of the form
\begin{equation}
\min_{u:\Omega\to\IM} F(u)
\end{equation}
with $\Omega \subset \R^d$ open and bounded. In this chapter, we will be concerned with problems
where  the image $u$ takes values in an $s$-dimensional \emph{manifold} $\IM$. 
Problems of this form are
wide-spread in image processing and especially in the processing of
manifold-valued images such as
InSAR \cite{massonnet1998_vogt},
EBSD \cite{bachmann2011_vogt},
DTI \cite{basser1994_vogt},
orientational/positional \cite{rosman2011_vogt}
data or images with values in non-flat color spaces such as hue-saturation-value
(HSV) or chromaticity-brightness (CB) color spaces \cite{chan2001_vogt}.

They come with
an inherent non-convexity, as the space of images $u\colon \Omega \to \IM$ is generally non-convex, with few exceptions, such as if $\IM$ is a Euclidean space, or if $\IM$ is a Hadamard manifold, if one allows for the more general notion of geodesic convexity \cite{bacak2014_vogt,bacak2016_vogt}. Except for these special cases, efficient and robust convex numerical optimization algorithms therefore cannot be applied and global optimization is generally out of reach.

\begin{figure}[t]
\tikzset{circle node/.style={circle,fill=black,inner sep=0,minimum size=#1}}
\begin{tikzpicture}[scale=0.43]
    \newlength{\crvogt} 
    \setlength{\crvogt}{2.7cm}
    \node at (-\crvogt,-1.35\crvogt) {};
    \node at (0,0) {$\IS^1$};
    \draw (0,0) circle (\crvogt);
    \node[circle node=5pt,label=53:{$x_1$}] at (53:\crvogt) {};
    \node[circle node=5pt,label=0:{$x_2$}] at (5:\crvogt) {};
    \node[circle node=3pt,label=29:{$\bar{x}$}] at (29:\crvogt) {};
    \node[circle node=3pt,label=-151:{$y$}] at (-151:\crvogt) {};
\end{tikzpicture}
\hfill\hspace{-0.5cm}
\begin{tikzpicture}[pin distance=13pt,every pin edge/.style={<-,>=latex},]
    \node at (-0.1,0.05)
        {\includegraphics[width=0.29\textwidth]{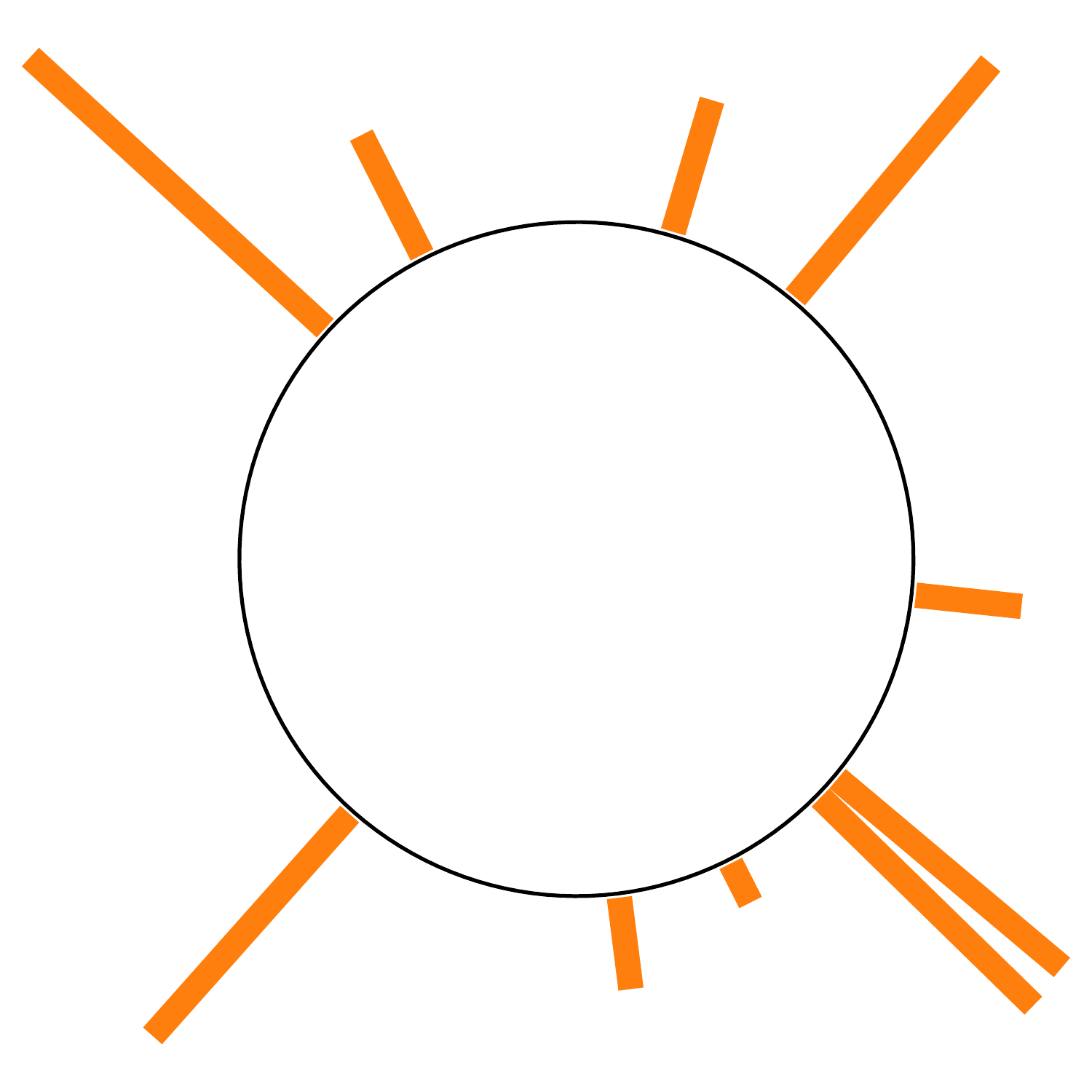}};
    \node[circle node=3pt,inner sep=0,pin={0:{$x_\text{start}$}}]
        at ({deg(pi)}:1.05) {};
    \node[circle node=3pt,inner sep=0,pin={-20:{$x_\text{local}$}}]
        at ({deg(3.4834278686879987)}:1.05) {};
    \node[circle node=3pt,inner sep=0,pin={160:{$\bar{x}$}}]
        at ({deg(0.22046264235717736)}:1.05) {};
\end{tikzpicture}
\hfill
\begin{tikzpicture}[pin distance=13pt,every pin edge/.style={<-,>=latex},]
    \node at (1.55,1.23)
        {\includegraphics[width=0.34\textwidth]{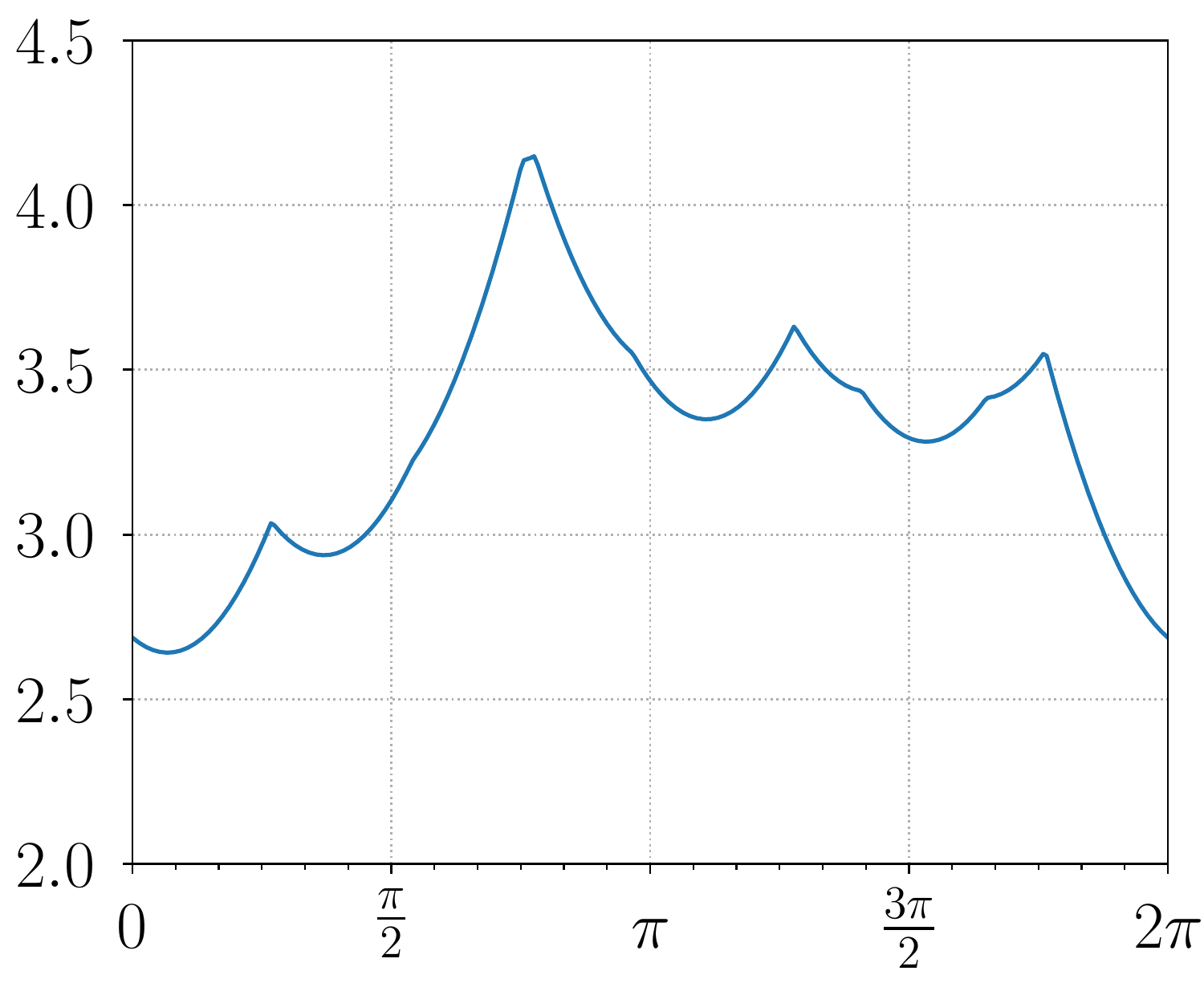}};
    \begin{scope}[scale=0.546]
    \node[inner sep=0,pin={95:{$x_\text{start}$}}] at (pi,0) {};
    \node[inner sep=0,pin={85:{$x_\text{local}$}}] at (3.4834278686879987,0) {};
    \node[inner sep=0,pin={85:{$\bar{x}$}}] at (0.22046264235717736,0) {};
    \end{scope}
\end{tikzpicture}
\hfill
\caption{%
    Variational problems where the feasible set is a non-Euclidean manifold are prone to local minima and non-uniqueness, which makes them generally much harder than their counterparts in $\R^n$. The example shows the generalization of the (weighted) mean to manifolds: the Riemannian center of mass $\bar{x}$ of points $x_i$ on a manifold -- in this case, the unit circle  $\IS^1$ -- is defined as the minimizer (if it exists and is unique) of the problem $\inf_{x\in\IS^1} \sum_{i} \lambda_id(x_i,x)^2$, where $d$ is the geodesic (angular) distance and $\lambda_i>0$ are given weights.
    \textbf{Left:} Given the two points $x_1$ and $x_2$, the energy for computing their ``average'' has a local minimum at $y$ in addition to the global minimum at $\bar{x}$. Compare this to the corresponding problem in $\R^n$, which has a strictly convex energy with the unique and explicit solution $(x_1+x_2)/2$.
    \textbf{Center and right:} When the number of points is increased and non-uniform weights are used (represented by the locations and heights of the orange bars), the energy structure becomes even less predictable. The objective function (right, parametrized by angle) exhibits a number of non-trivial local minimizers that are not easily explained by global symmetries. Again, the corresponding problem -- computing a weighted mean -- is trivial in $\R^n$.
    Starting from $x_\text{start} = \pi$, our functional lifting implementation finds the global minimizer~$\bar{x}$, while gradient
    descent (a local method) gets stuck in the local minimizer $x_\text{local}$.
    Empirically, this behaviour can be observed for any other choice of points
    and weights, but there is no theoretical result in this direction.
%
%
%
%
}\label{fig:circle-riem-com_vogt}
\end{figure}
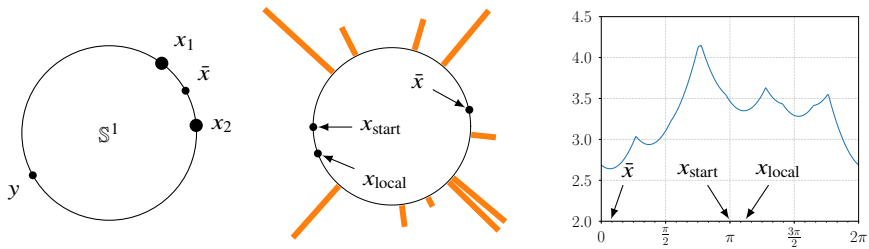

The inherent non-convexity of the feasible set is not only an issue of representation. Even for seemingly simple problems, such as the problem of computing the  Riemannian center of mass for a number of points on the unit circle, it can affect the energy in surprisingly intricate ways, creating multiple local minimizers and non-uniqueness 
(Fig.~\ref{fig:circle-riem-com_vogt}).
The equivalent operation in Euclidean space, computing the  weighted mean, is a simple convex (even linear)\ operation, with a unique, explicit solution. 

The problem of non-convexity is not unique to our setting, but rather ubiquitous in a much broader context of image and signal processing:
amongst others, image segmentation, 3D reconstruction, image matching,
optical flow and image registration, superresolution, inpainting, edge-preserving image restoration with the Mumford-Shah and Potts model, machine learning, and many statistically or physically motivated models involve intrinsically non-convex feasible sets or energies.
When applied to such non-convex problems, local optimization strategies often get stuck in local minimizers.

In \emph{convex relaxation} approaches, an energy functional is approximated by a
convex one whose global optimum can be found numerically and whose minimizers lie within a small neighborhood around the actual
solution of the problem. A popular convex relaxation technique that applies to a wide range of problems
from image and signal processing is \emph{functional lifting}.
With this technique, the feasible set is embedded into a higher-dimensional space
where efficient convex approximations of the energy functional are easier
available.

\textbf{Overview and contribution.}
In the following sections, we will give a brief introduction to the concept of functional lifting and explore its generalization to manifold-valued problems. Our aim is to provide a survey-style introduction to the area, therefore we will provide references and numerical experiments on the way. In contrast to prior work, we will explain existing results in an updated finite element-based framework. Moreover, we propose extensions to handle general regularizers other than the total variation on manifolds, and to apply the ``sublabel-accurate'' methods to manifold-valued problems.

\subsection{Functional lifting in Euclidean spaces}
The problem of finding a function $u\colon \Omega \to \Gamma$ that assigns
a \emph{label} $u(x) \in \Gamma$ from a 
\emph{discrete} range
$\Gamma$ to each point $x$ in a continuous domain $\Omega \subset \R^d$,
while minimizing an energy function $F(u)$, is commonly called a \emph{continuous
multi-label (or multi-class labeling) problem} in the image processing community
\cite{pock2008_vogt,lellmann2009_vogt}.
The name comes from the interpretation of this setting as the continuous
counterpart to the fully discrete problem of assigning to each vertex of a graph
one of finitely many labels $\gamma_1,\dots,\gamma_L$ while minimizing a given
cost function \cite{greig1989_vogt,calinescu1998_vogt,kleinberg2002_vogt,ishikawa2003_vogt}.

The prototypical application of multi-labeling techniques is multi-class image
segmentation, where the task is to partition a given image into finitely many
regions.
In this case, the label set $\Gamma$ is discrete and each label represents one
of the regions so that $u^{-1}(\gamma) \subset \Omega$ is the region that
is assigned label $\gamma$.

In the fully discrete setting, one way of tackling first-order multi-label
problems is to look for good linear programming relaxations
\cite{calinescu1998_vogt,kleinberg2002_vogt,ishikawa2003_vogt}. These approaches were subsequently translated to continuous domains $\Omega$ for the two-class \cite{chan2006_vogt}, multi-class \cite{zach2008_vogt,pock2008_vogt,lellmann2011c_vogt,bae2011_vogt}, and vectorial \cite{goldluecke2010_vogt} case, resulting in non-linear, but convex, relaxations. By honoring the continuous nature of $\Omega$, they reduce metrication errors and improve isotropy \cite{strekalovskiy2011_vogt,strekalovskiy2012_vogt,goldluecke2013_vogt,strekalovskiy2015_vogt}, see \cite{lellmann2013b_vogt} for a discussion and more references.

The general strategy, which we will also follow for the manifold-valued case, is to replace the energy minimization problem
\begin{equation}
    \min_{u\colon \Omega \to \Gamma } F(u),\label{eq:intro-unlifted_vogt}
\end{equation}
by a problem
\begin{equation}
    \min_{v\colon \Omega \to X} \tilde F(v),\label{eq:intro-lifted_vogt}
\end{equation}
where $X$ is some ``nice'' convex set of larger dimension than $\Gamma$ with the property
that there is an embedding $i\colon \Gamma \hookrightarrow X$ and
$F(u) \approx \tilde F(i \circ u)$ in some sense whenever
$u\colon \Omega \to \Gamma$.

In general, the lifted functional $\tilde F$ is chosen in such a way that it exhibits
favorable (numerical or qualitative) properties compared with the original
functional $F$ while being sufficiently close to the original functional so that
minimizers of $\tilde F$ can be expected to have some recoverable relationship
with global minimizers of $F$.
Usually, $\tilde F$ is chosen to be convex when $F$ is not, which will make
the problem amenable for convex optimization algorithms and allows to find a global minimizer of the lifted problem.

While current lifting strategies generally avoid local minimizers of the original problem, they are still an approximation and they are generally not guaranteed to find the global minimizers of the original problem.

A central difficulty is that some simplifications have to be performed in the lifting process in order to make it computationally feasible, which may lose information about the original problem. As a result,
global minimizers $v\colon \Omega \to X$  of the lifted problem need not be
in the image of $\Gamma$ under the embedding $i\colon \Gamma \hookrightarrow X$ and therefore are not directly associated with a function in the original space.

The process of projecting a solution back to the original space of functions
$u\colon \Omega \to \Gamma$ is a difficult problem and, unless $\Gamma$ is
scalar \cite{pock2010_vogt}, the projection cannot be expected to be a minimizer
of the original functional (see the considerations in
\cite{federer1974_vogt,lavenant2017_vogt,vogt2019_vogt}).
These difficulties may be related to the fact that the original problems are
NP-hard \cite{cremers2012_vogt}.
As in the discrete labeling setting \cite{kleinberg2002_vogt}, so-called
rounding strategies have been proposed in the continuous case
\cite{lellmann2012_vogt,lellmann2011phd_vogt} that come with an \emph{a priori} bound for
the relative gap between the minimum of the original functional and the value
attained at the projected version of a minimizer to the lifted functional. For the manifold-valued case considered here, we are not aware of a similar result yet.

In addition to the case of a \emph{discrete} range $\Gamma$, relaxation methods have been derived for dealing with a \emph{continuous} (non-discrete) range, most notably the scalar case $\Gamma\subseteq\R$ \cite{alberti2003_vogt,pock2010_vogt}. They typically consider \emph{first-order} energies that depend
pointwise on $u$ and $\nabla u$ only:
\begin{equation}\label{eq:first-order-cml_vogt}
    F(u) = \int_\Omega f(x,u(x),\nabla u(x))\,dx.
\end{equation}
The equivalent problem class in the fully discrete setting consists of the energies with only unary (depending on one vertex's label) and pairwise (depending on two vertices' labels) terms.

For the problem \eqref{eq:first-order-cml_vogt}, applying a strategy as in \eqref{eq:intro-unlifted_vogt}--\eqref{eq:intro-lifted_vogt} comes with a substantial increase in dimensions. These relaxation approaches therefore have been called \emph{functional lifting},
starting from the paper \cite{pock2009_vogt} where the (non-convex) Mumford-Shah
functional for edge-preserving image regularization and segmentation is lifted
to a space of functions $v\colon \Omega \times \Gamma \to [0,1]$,
$\Gamma \subset \R$.
The authors use the special ``step function'' lifting $X = \{ v\colon \Gamma \to [0,1] \}$ and
$i(z^*) = v$ with $v(z) = 1$ if $z \leq z^*$ and $0$ otherwise, which is only available in the scalar case

In this case, the integrand
$f\colon \Omega \times \Gamma \times \R^{s,d} \to \R$ in
\eqref{eq:first-order-cml_vogt} is assumed to be convex in the third component
and nonnegative. The less restrictive property of polyconvexity has been
shown to be sufficient \cite{windheuser2016_vogt,mollenhoff2019_vogt}, so that also minimal surface problems fit into this framework.
The continuous formulations can be demonstrated \cite{pock2009_vogt,mollenhoff2019_vogt}
to have strong connections with the method of calibrations \cite{alberti2003_vogt}
and with the theory of currents \cite{giaquinta1998_vogt}. 


In this paper, we will consider the more general case of $\Gamma = \IM$ having a manifold structure. We will also restrict ourselves to first-order models.
Only very recently, attempts at generalizing the continuous lifting
strategies to models with higher-order regularization have been made -- for
regularizers that depend on the Laplacian \cite{loewenhauser2018_vogt,vogt2019_vogt}
in case of vectorial ranges $\Gamma \subset \R^s$ and for the total generalized
variation \cite{ranftl2013_vogt,strecke2019_vogt} in case of a scalar range $\Gamma \subset \R$.
However, in contrast to the first-order theory, the higher-order models, although empirically useful, are still considerably less mathematically validated.
Furthermore, we mention that there are models where the image \emph{domain}
$\Omega$ is replaced by a shape (or manifold) \cite{delaunoy2009_vogt,bernard2017_vogt}, which is beyond the scope of this survey.

\subsection{Manifold-valued functional lifting}

\begin{figure}[t]
\definecolor{lightblue}{RGB}{180,220,255}
\begin{tikzpicture}[scale=0.42]
    \shade[left color=lightblue!40!white,right color=lightblue]
      (0,0) to[out=-10,in=150] (6,-2) -- (12,1) to[out=150,in=-10] (5.5,3.7) -- cycle;
      
    \node[circle,fill=red,inner sep=0,minimum size=3pt,label=-153:{$z$}] at (5.2,1.0) {};
      
    \draw[->] (11.5,2) to[out=20,in=160] node[above] {$i\colon \IM \to \IP(\IM)$} (17.3,2);
    
    \begin{scope}[shift={(16,0)}]
    \shade[left color=lightblue!40!white,right color=lightblue]
      (0,0) to[out=-10,in=150] (6,-2) -- (12,1) to[out=150,in=-10] (5.5,3.7) -- cycle;
      
    \draw[->,thick,red,>=latex] (5.2,1.0) -- ++(0,4.5);
    \shade[left color=lightblue!30!black,right color=lightblue!90!white]
        (5.2,1.05) -- ++(0,-0.1) -- ++(4.95,1.0) -- ++(0,0.1) -- cycle;
    \node[circle,fill=red,inner sep=0,minimum size=2pt,label=-90:{$i(z) = \delta_z$}] at (5.2,1.0) {};
    \end{scope}
\end{tikzpicture}
\caption{%
    A manifold $\IM$ is embedded into the space $\IP(\IM)$ of probability
    measures via the identification of a point $z \in \IM$ with the Dirac
    point measure $\delta_z$ concentrated at $z$. This ``lifts'' the problem into a higher-dimensional \emph{linear} space, which is much more amenable to global optimization methods. 
}\label{fig:prob-embedding_vogt}
\end{figure}
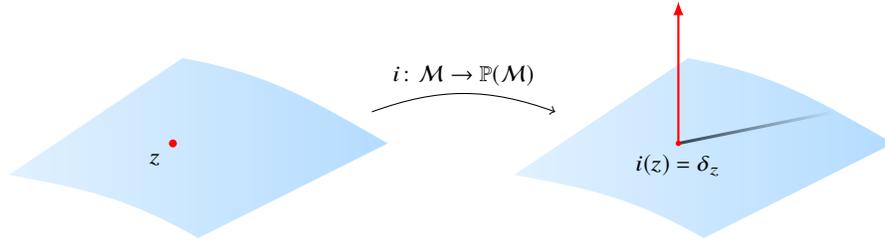

In this chapter, we will be concerned with problems where $\Gamma$ has a manifold structure.
The first step towards applying lifting methods to such problems was an application to the restoration of cyclic data \cite{strekalovskiy2011b_vogt,cremers2012_vogt} with $\Gamma = \IS^1$,
which was later \cite{lellmann2013_vogt} generalized for
the case of total variation regularization to data with values in more general
manifolds.
In \cite{lellmann2013_vogt}, the functional lifting approach is applied to a
first-order model with  total variation regularizer,
\begin{equation}
    F(u) = \int_\Omega \rho(x,u(x)) dx + \lambda \TV(u),
\end{equation}
for $u\colon \Omega \to \IM$, where $\Gamma = \IM$ is an $s$-dimensional
manifold and $\rho\colon \Omega \times \IM \to \R$ is a pointwise data
discrepancy.
The lifted space is chosen to be $X = \IP(\IM)$, the space of Borel probability
measures over $\IM$, with embedding $i\colon \IM \hookrightarrow \IP(\IM)$, where
$i(z) := \delta_z$ is the Dirac point measure with unit mass concentrated at
$z \in \IM$ (see Fig.~\ref{fig:prob-embedding_vogt}).
The lifted functional is
\begin{equation}
    \tilde F(v) = \int_\Omega \langle \rho(x,\cdot), v(x) \rangle \,dx
        + \lambda \widetilde{\TV}(v),
\end{equation}
where $\langle g, \mu \rangle := \int_\IM g \,d\mu$ for $g \in C(\IM)$
and $\mu \in \IP(\IM)$.
Furthermore,
\begin{equation}
    \widetilde{\TV}(v) := \sup\left\{
        \int_{\Omega} \langle \Div_x p(x,\cdot), v(x) \rangle \,dx :
        p\colon \Omega \times \IM \to \R,
        \|\nabla_z p\|_{\infty} \leq 1
    \right\}.
\end{equation}
The Lipschitz constraint $\|\nabla_z p\|_{\infty} \leq 1$, where
\begin{equation}
    \|\nabla_z p\|_\infty := \sup \left\{\|\nabla_z p(x,z)\|_{\sigma,\infty}:
        (x,z) \in \Omega \times \IM
    \right\},
\end{equation}
and $\|\cdot\|_{\sigma,\infty}$ the spectral (operator) norm, can be explained by
a functional analytic perspective \cite{vogt2018_vogt} on this lifting strategy:
The lifted total variation functional is the vectorial total variation semi-norm
for functions over $\Omega$ with values in a certain Banach space of measures.
The topological dual space of this space of measures is the space of Lipschitz
continuous functions over $\IM$.
However, this interpretation does not generalize easily to other regularizers.
We will instead base our model for general convex regularizers on the theory
of currents as presented in \cite{mollenhoff2019_vogt}.

\textbf{Sublabel accuracy}.
\begin{figure}[t]
\begin{tabular}{cccc}
\includegraphics[width=0.24\textwidth,trim=189 100 165 100,clip]{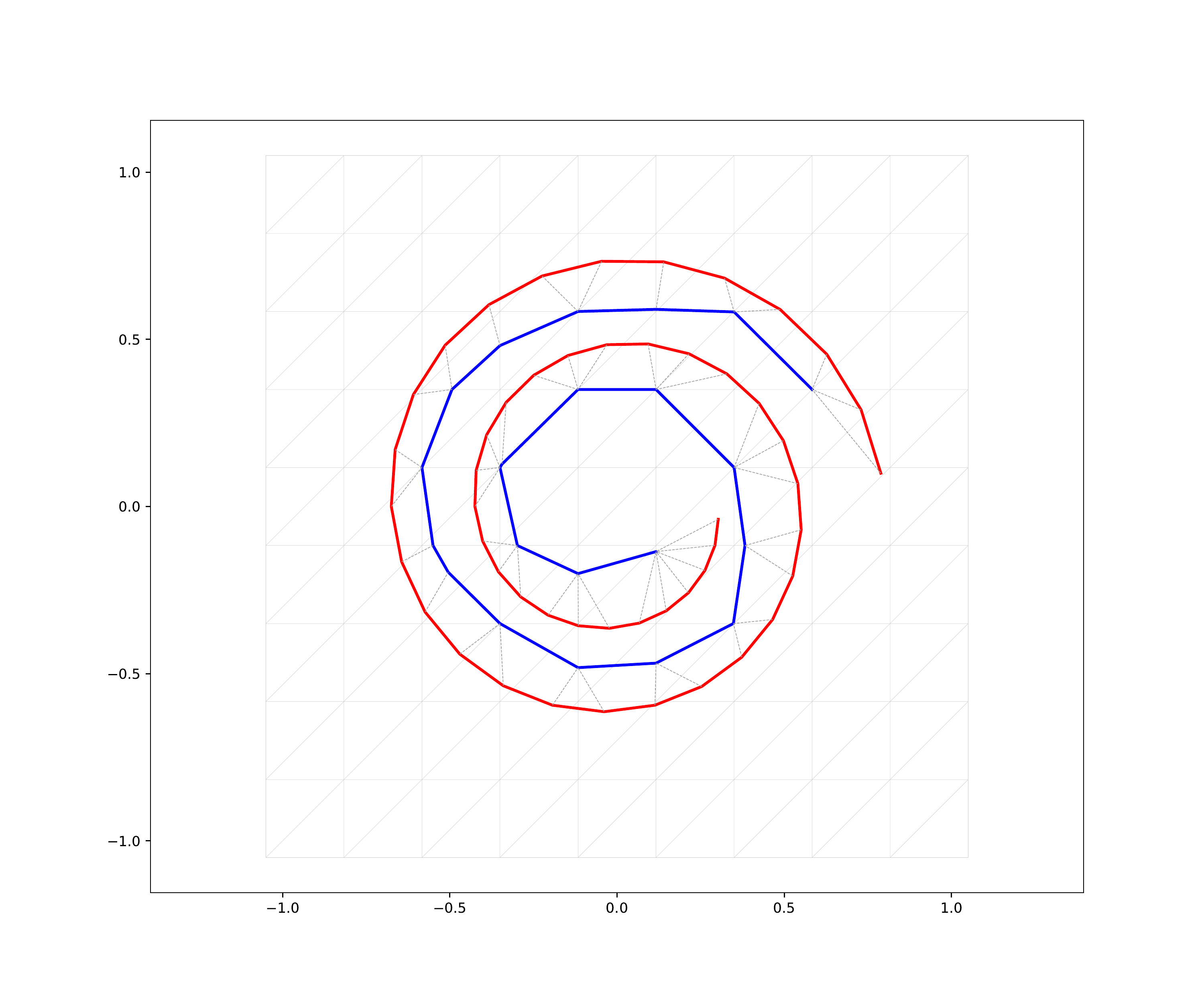}&%
\includegraphics[width=0.24\textwidth,trim=189 100 165 100,clip]{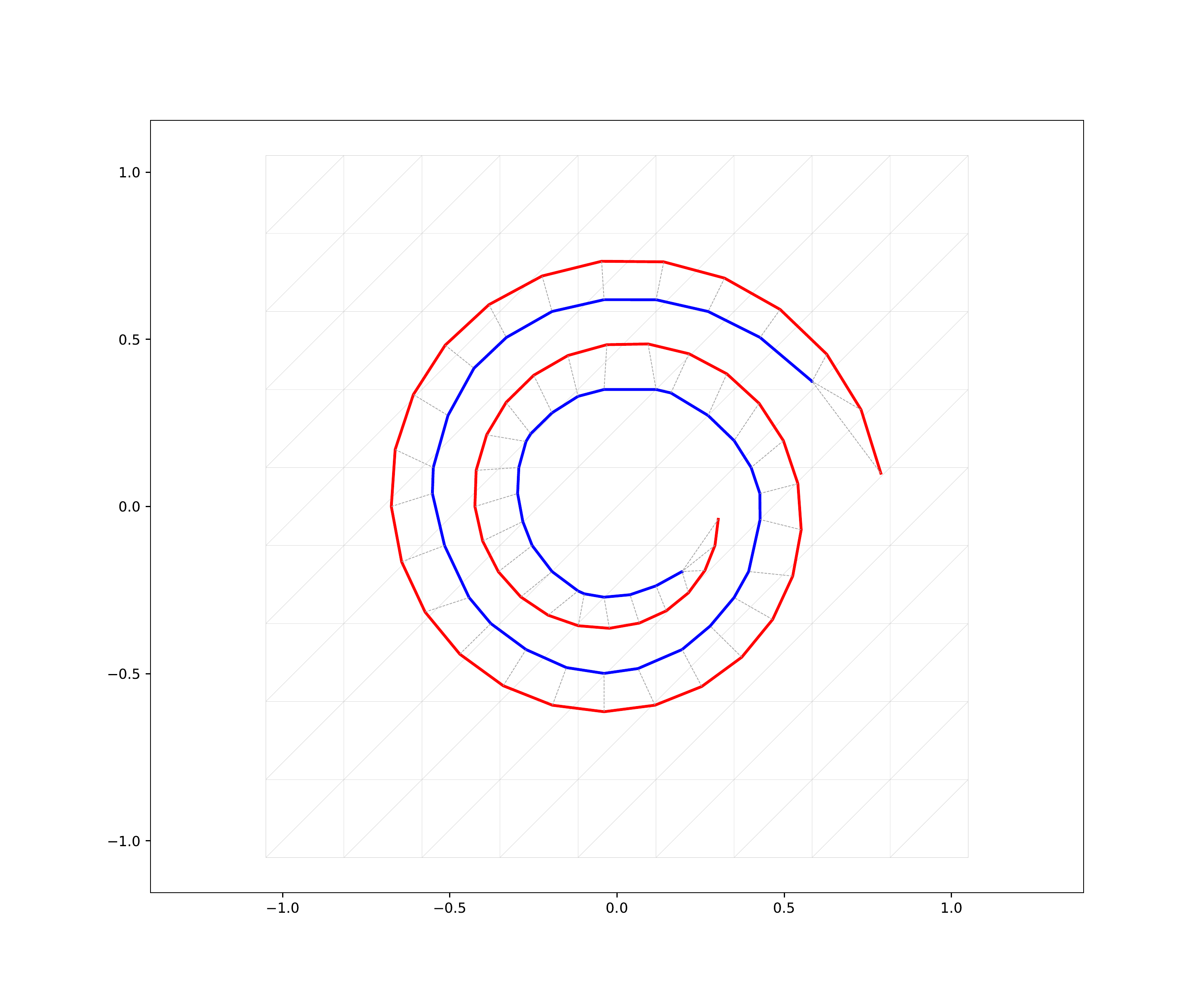}&%
\includegraphics[width=0.24\textwidth,trim=189 100 165 100,clip]{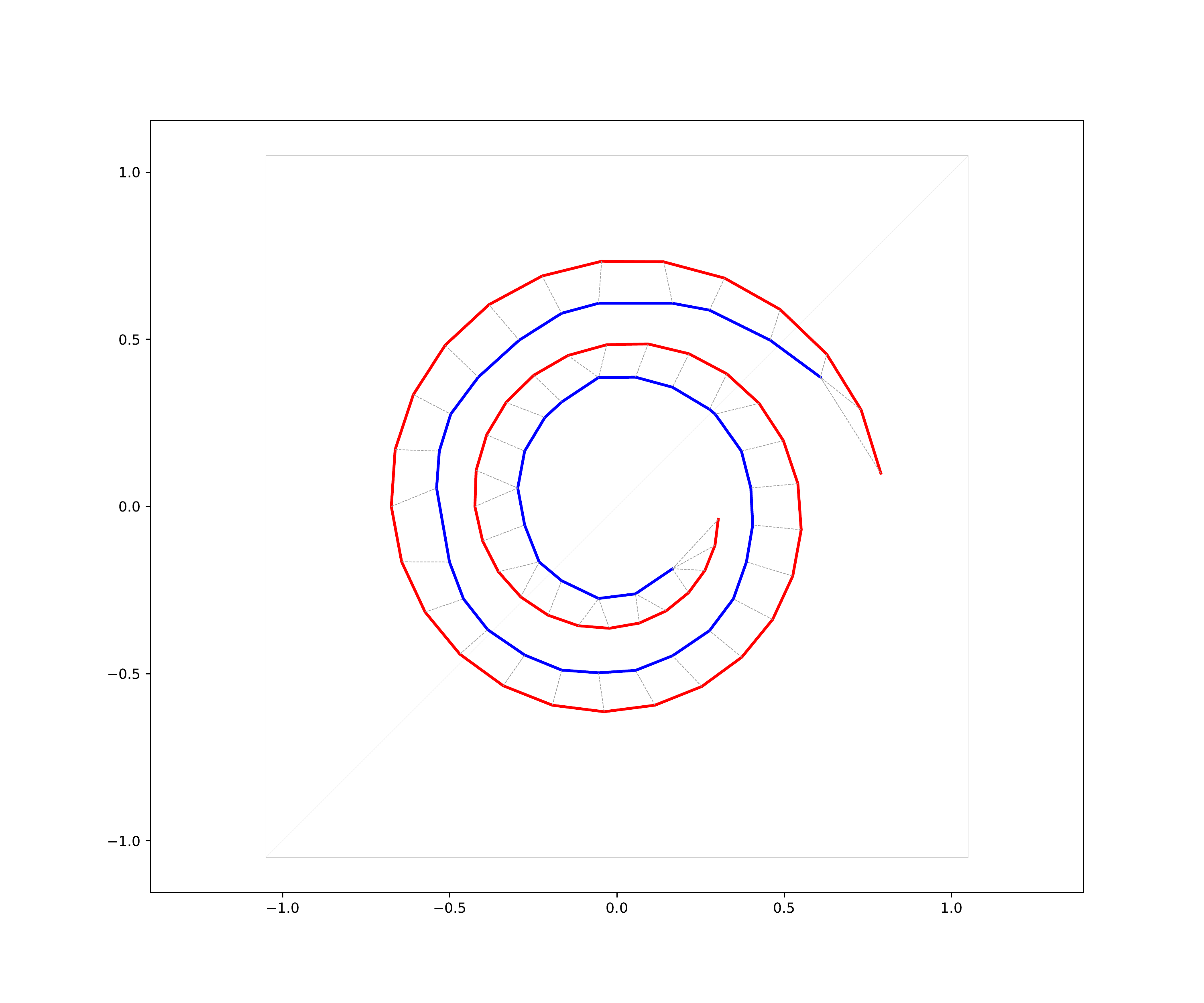}&%
\includegraphics[width=0.24\textwidth,trim=189 100 165 100,clip]{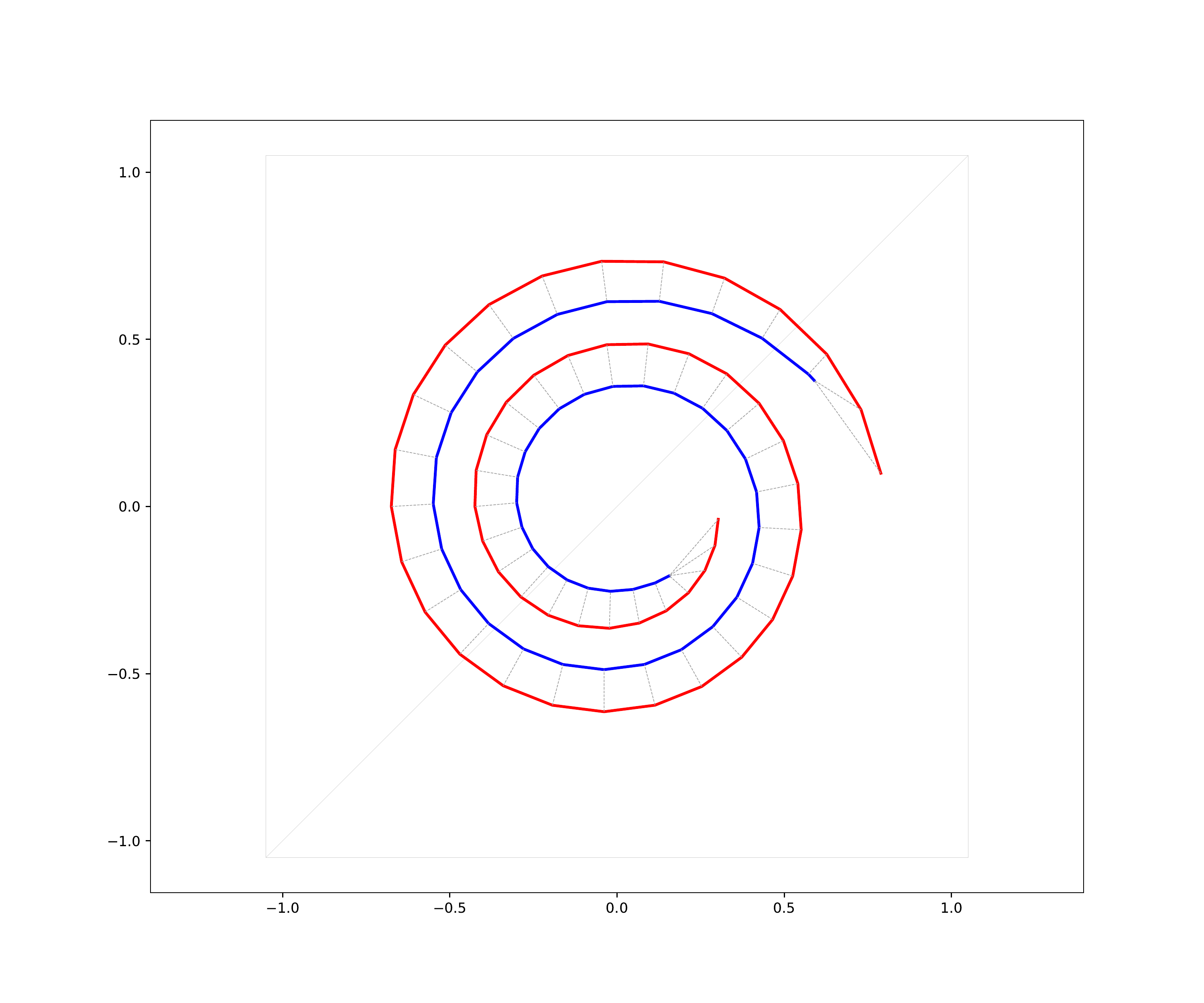}\\
\includegraphics[width=0.24\textwidth,trim=358 384 280 167,clip]{fig/plot-flat-1d-10-2.pdf}&%
\includegraphics[width=0.24\textwidth,trim=358 384 280 167,clip]{fig/plot-flat-1d-10-4.pdf}&%
\includegraphics[width=0.24\textwidth,trim=358 384 280 167,clip]{fig/plot-flat-1d-2-20.pdf}&%
\includegraphics[width=0.24\textwidth,trim=358 384 280 167,clip]{fig/plot-flat-1d-2-exact.pdf}\\
\cite{lellmann2013_vogt}, $10\times 10$ labels&
\cite{laude2016_vogt}, $10\times 10$ labels&
\cite{laude2016_vogt}, $2\times 2$ labels&
\cite{laude2016_vogt}, $2\times 2$ labels\\
label bias&
no label bias&
sublabel-accurate&
exact data term
\end{tabular}
\caption{%
    Rudin-Osher-Fatemi (ROF) $L^2-\TV$ denoising (blue) of an (Euclidean) vector-valued signal
    $u\colon [0,1] \to \R^2$ (red), visualized as a curve in the flat manifold $\IM=\R^2$.
    The problem is solved by the continuous multi-labeling framework with
    functional lifting described in this chapter. 
    The discretization points (labels) in the range $\IM$, which are necessary for the implementation
    of the lifted problem, are visualized by the gray grid. \textbf{Left:}
    The method proposed in \cite{lellmann2013_vogt} does not force the solution
    to assume values at the grid points (labels), but still shows significant bias towards edges of the grid (blue curve). \textbf{Second from left:}
    With the same number of labels, the method from \cite{laude2016_vogt} is
    able to reduce label bias by improving data term discretization. \textbf{Second from right:}
    Furthermore, the method from \cite{laude2016_vogt} allows to exploit the convexity of the data term to get decent
    results with as little as four grid points.
    \textbf{Right:} Further exploiting the quadratic form of the data term even
    produces the numerically exact reference solution, which in this case can be precisely computed using the unlifted formulation due to the convexity of the problem. This shows that for the Euclidean fully convex case, the sublabel-accurate lifting allows to recover the exact solution with careful discretization.
}\label{fig:label-bias-flat_vogt}
\end{figure}%
\begin{figure}[t]
%
\includegraphics[width=0.33\textwidth]{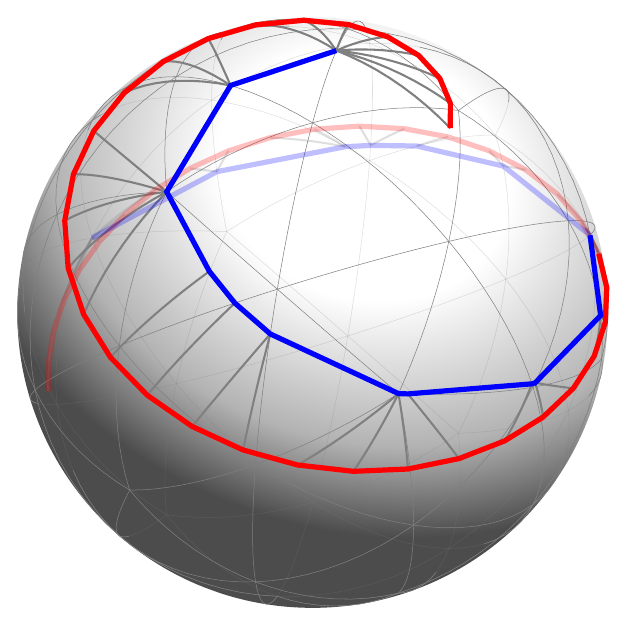}
\includegraphics[width=0.33\textwidth]{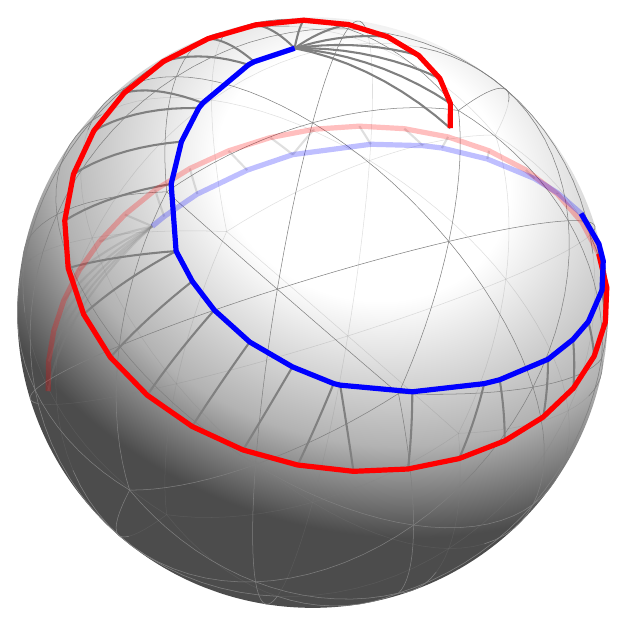}
\includegraphics[width=0.33\textwidth]{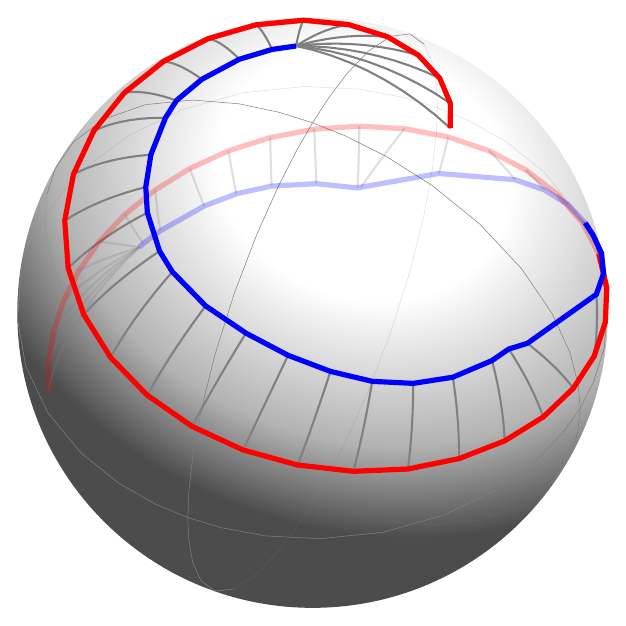}
\caption{%
    Total Variation denoising (blue) of a signal
    $u\colon [0,1] \to \IS^2$ with values in $\IS^2$ (red), visualized as curves
    on the two-dimensional sphere embedded into $\R^3$.
    The problem is solved by the continuous multi-labeling framework with
    functional lifting described in this chapter. 
    The discretization points (labels), that are necessary for the implementation
    of the lifted problem, are visualized by the gray grid.
    \textbf{Left:} The method proposed in \cite{lellmann2013_vogt} does not force the solution
    to take values at the grid points, but still shows significant grid bias.
    \textbf{Center:} With the same number of labels, our proposed method, motivated by
    \cite{laude2016_vogt}, reduces label bias by improving data term
    discretization.
    \textbf{Right:} Furthermore, our method can get excellent results with as little as 6 grid points
    (right).
    Note that the typical contrast reduction that occurs in the classical Euclidean ROF can also be observed in the manifold-valued case in the form of a shrinkage towards the Fr\'{e}chet mean.
}\label{fig:label-bias-s2_vogt}
\end{figure}
While the above model comes with a fully continuous description, a numerical
implementation requires the discretization of $\Omega$ as well as the range $\Gamma$.
This introduces two possible causes for errors: \emph{metrication errors} and \emph{label bias}.

\emph{Metrication errors} are artifacts related to the graph or grid representation of the spatial image domain $\Omega$, finite difference operators, and the choice of metric thereof. They manifest mostly in unwanted anisotropy, missing rotational invariance, or blocky diagonals. They constitute a common difficulty with all variational problems and lifting approaches \cite{klodt2008_vogt}.

In contrast, \emph{label bias} means that the discretization favors solutions that assume values at the chosen ``labels'' (discretization points) $Z^1,\dots,Z^L$ in the \emph{range} $\Gamma$ (see Fig.~\ref{fig:label-bias-flat_vogt} and~\ref{fig:label-bias-s2_vogt}).
This is very desirable for discrete $\Gamma$, but in the context of manifolds, severely limits accuracy and forces a suitably fine discretization of the range.


In more recent so-called \emph{sublabel-accurate} approaches for scalar and
vectorial ranges $\Gamma$, more emphasis is put on the discretization
\cite{zach2012_vogt,mollenhoff2016_vogt,laude2016_vogt} to get rid of label bias
in models with total variation regularization, which allows to greatly reduce the number of discretizations points for the range $\Gamma$.
In a recent publication \cite{mollenhoff2017_vogt}, the gain in sublabel accuracy
is explained to be caused by an implicit application of first-order finite
elements on $\Gamma$ as opposed to previous approaches that can be interpreted as using zero-order elements,
which naturally introduces label-bias.
An extension of the sublabel-accurate approaches to arbitrary convex regularizers
using the theory of currents was recently proposed in \cite{mollenhoff2019_vogt}.

Motivated by these recent advances, we propose to extend the methods from
\cite{lellmann2013_vogt} for manifold-valued images to arbitrary convex
regularizers, making use of finite element techniques on manifolds
\cite{dziuk2013_vogt}.
This reduces label bias and thus the amount of labels necessary in the
discretization.

\subsection{Further related work}

The methods proposed in this work are applicable to variational problems with
values in manifolds of dimension $s \leq 3$.
The theoretical framework applies to manifolds of arbitrary dimension, but the
numerical costs increase exponentially for dimensions $4$ and larger.

An alternative is to use \emph{local} optimization methods on manifolds. A reference for the smooth case is \cite{absil2009_vogt}. For non-smooth energies, methods such as the cyclic proximal point, Douglas-Rachford, ADMM and (sub-)gradient descent algorithm have been applied to first and second
order TV and TGV as well as Mumford-Shah and Potts regularization approaches in
\cite{%
    weinmann2014_vogt,%
    weinmann2015_vogt,%
    baust2016_vogt,%
    bergmann2016_vogt,%
    bredies2018_vogt,%
    bergmann2018b_vogt%
}.
These methods are generally applicable to manifolds of any
dimension whose (inverse) exponential mapping can be evaluated in reasonable
time and quite efficient in finding a local miminum, but can get stuck in local extrema.
Furthermore, the use of total variation regularization in these frameworks is currently limited to anisotropic formulations; Tikhonov regularization
was proposed instead for isotropic regularization
\cite{weinmann2014_vogt,bergmann2018c_vogt}.
An overview of applications, variational models and local optimization methods
is given in \cite{bergmann2018c_vogt}.

Furthermore, we mention that, beyond variational models, there exist statistical
\cite{fletcher2012_vogt}, discrete graph-based \cite{bergmann2018_vogt},
wavelet-based \cite{storath2018_vogt}, PDE-based \cite{chefdhotel2004_vogt} and
patch-based \cite{laus2017_vogt} models for the processing and regularization of
manifold-valued signals.

%% file: chpt/sec1.tex
\section{Submanifolds of $\R^N$}

We formulate our model for submanifolds of $\R^N$ which is no restriction
by the Whitney embedding theorem \cite[Thm.~6.15]{lee2013_vogt}.
For an $s$-dimensional submanifold of $\R^N$ and $\Omega \subset \R^d$
open and bounded, differentiable functions $u\colon\Omega \to \IM$ are regarded
as a subset of differentiable functions with values in $\R^N$.
For those functions, a Jacobian $Du(x) \in \R^{N,d}$ in the Euclidean sense
exists that can be identified with the push-forward of the tangent space
$T_x\Omega$ to $T_{u(x)}\IM$, i.e., for each $x \in \Omega$ and
$\xi \in \R^d = T_x\Omega$, we have
\begin{equation}
    Du(x)\xi \in T_{u(x)}\IM \subset T_{u(x)}\R^N.
\end{equation}
On the other hand, for differentiable maps $p\colon \IM \to \R^d$, there exists
an extension of $p$ to a neighborhood of $\IM \subset \R^N$ that is constant in
normal directions and we denote by $\nabla p(z) \in \R^{N,d}$ the Jacobian of
this extension evaluated at $z \in \IM$.
Since the extension is assumed to be constant in normal directions, i.e.,
$\nabla p(z)\zeta = 0$ whenever $\zeta \in N_z \IM$ (the orthogonal
complement of $T_z \IM$ in $\R^N$), this definition is independent of the choice
of extension.

\subsection{Calculus of Variations on submanifolds}

In this section, we generalize the total variation based approach in
\cite{lellmann2013_vogt} to less restrictive first-order variational problems
by applying the ideas from functional lifting of vectorial problems
\cite{mollenhoff2019_vogt} to manifold-valued problems. Most derivations will be formal; we leave a rigorous choice of function spaces as well as an analysis of well-posedness for future work.
We note that theoretical work is available for the scalar-valued case in
\cite{alberti2003_vogt,pock2010_vogt,bouchitte2018_vogt} and for the vectorial
and for selected manifold-valued cases in \cite{giaquinta1998_vogt}.

We consider variational models on functions $u\colon\Omega \to \IM$,
\begin{equation}
    F(u) := \int_\Omega f(x,u(x),Du(x)) \,dx,
\end{equation}
for which the integrand $f\colon \Omega \times \IM \times \R^{N,d} \to \R$
is convex in the last component.
Note that the dependence of $f$ on the full Jacobian of $u$ spares us dealing
with the tangent bundle push-forward $T\Omega \to T\IM$ in a coordinate-free
way, thus facilitating discretization later on.

Formally, the lifting strategy for vectorial problems proposed in
\cite{mollenhoff2019_vogt} can be generalized to this setting by replacing
the range $\Gamma$ with $\IM$.
As the lifted space, we consider the space of probability measures on the
Borel $\sigma$-Algebra over $\IM$, $X = \IP(\IM)$, with embedding
$i\colon \IM \to \IP(\IM)$, where $i(z) = \delta_z$ is the Dirac
point mass concentrated at $z \in \IM$.
Furthermore, we write $\Sigma := \Omega \times \IM$ and, for
$(x,z) = y \in \Sigma$, we define the coordinate projections $\pi_1 y := x$ and $\pi_2 y := z$.
Then, for $v\colon \Omega \to \IP(\IM)$, we define the lifted functional
\begin{equation}
    \tilde F(v) := \sup\left\{
        \int_{\Omega} \langle -\Div_x p(x,\cdot) + q(x,\cdot), v(x) \rangle \,dx :
        (\nabla_z p,q) \in \IK
    \right\},
\end{equation}
where $\langle g, \mu \rangle := \int_\IM g \,d\mu$ is the dual pairing between
$g \in C(\IM)$ and $\mu \in \IP(\IM)$ and
\begin{equation}
    \IK := \left\{
        (P,q) \in C(\Sigma; \R^{N,d} \times \R) :
        f^*(\pi_1 y,\pi_2 y,P(y)) + q(y) \leq 0 \,\forall y \in \Sigma
    \right\},
\end{equation}
where $f^*(x,z,\zeta) := \sup_\xi \langle \zeta, \xi \rangle - f(x,z,\xi)$ is
the convex conjugate of $f$ with respect to the last variable.

In the following, the integrand
$f\colon \Omega \times \IM \times \R^{N,d} \to \R$ is assumed to decompose as
\begin{equation}
    f(x,z,\xi) = \rho(x,z) + \eta(P_z\xi)
\end{equation}
into a pointwise data term $\rho\colon \Omega \times \IM \to \R$ and a convex
regularizer $\eta\colon \R^{s,d} \to \R$ that only depends on an $s$-dimensional
representation of vectors in $T_z \IM$ given by a surjective linear map
$P_z \in \R^{s,N}$ with $\kernel(P_z) = N_z \IM$.

This very general integrand covers most first-order
models in the literature on manifold-valued imaging problems. It applies in particular to isotropic and anisotropic regularizers that
depend on (matrix) norms of $Du(x)$ such as the Frobenius or spectral norm (or operator norm) where $P_z$ is taken to be an arbitrary orthogonal basis
transformation.
Since $z \mapsto P_z$ is not required to be continuous, it can also be applied
to non-orientable manifolds such as the Moebius strip or the Klein bottle where
no continuous orthogonal basis representation of the tangent bundle $T\IM$
exists.

Regularizers of this particular form depend on the manifold through
the choice of $P_z$ only.
This is important because we approximate $\IM$ in the course of our proposed
discretization by a discrete (simplicial) manifold $\IM_h$ and the tangent
spaces $T_z \IM$ are replaced by the linear spaces spanned by the simplicial
faces of $\IM_h$.

\subsection{Finite elements on submanifolds}

\begin{figure}[t]
    \centering
    
    \hfill
    \includegraphics[width=0.5\textwidth,clip,trim=20 130 20 50]%
        {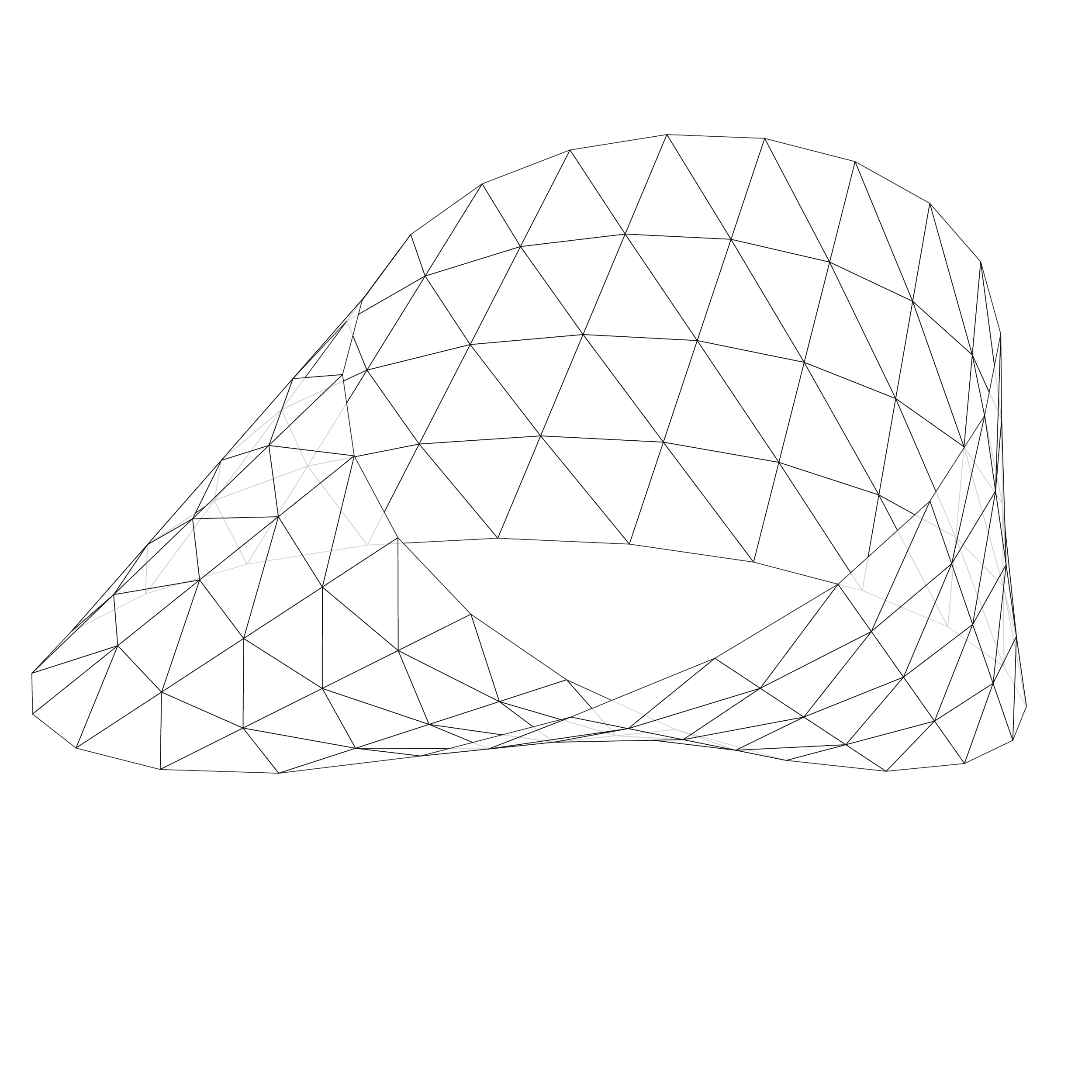}
    \hfill
    \includegraphics[width=0.4\textwidth]{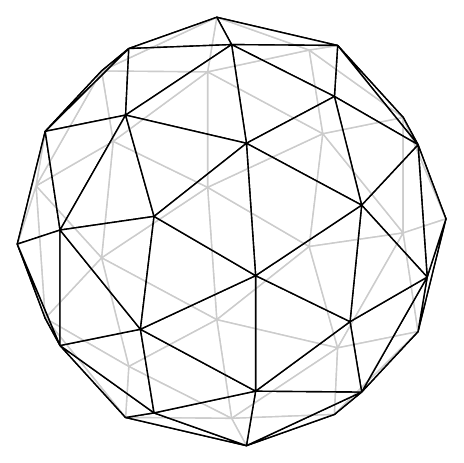}
    \hfill
    
    \caption{
        Triangulated approximations of the Moebius strip (left) and the
        two-dimensional sphere \t(right) as surfaces embedded into $\R^3$.
    }\label{fig:trisurf_vogt}
\end{figure}

\begin{figure}[t]
\sidecaption[t]
    \includegraphics[width=0.64\textwidth]{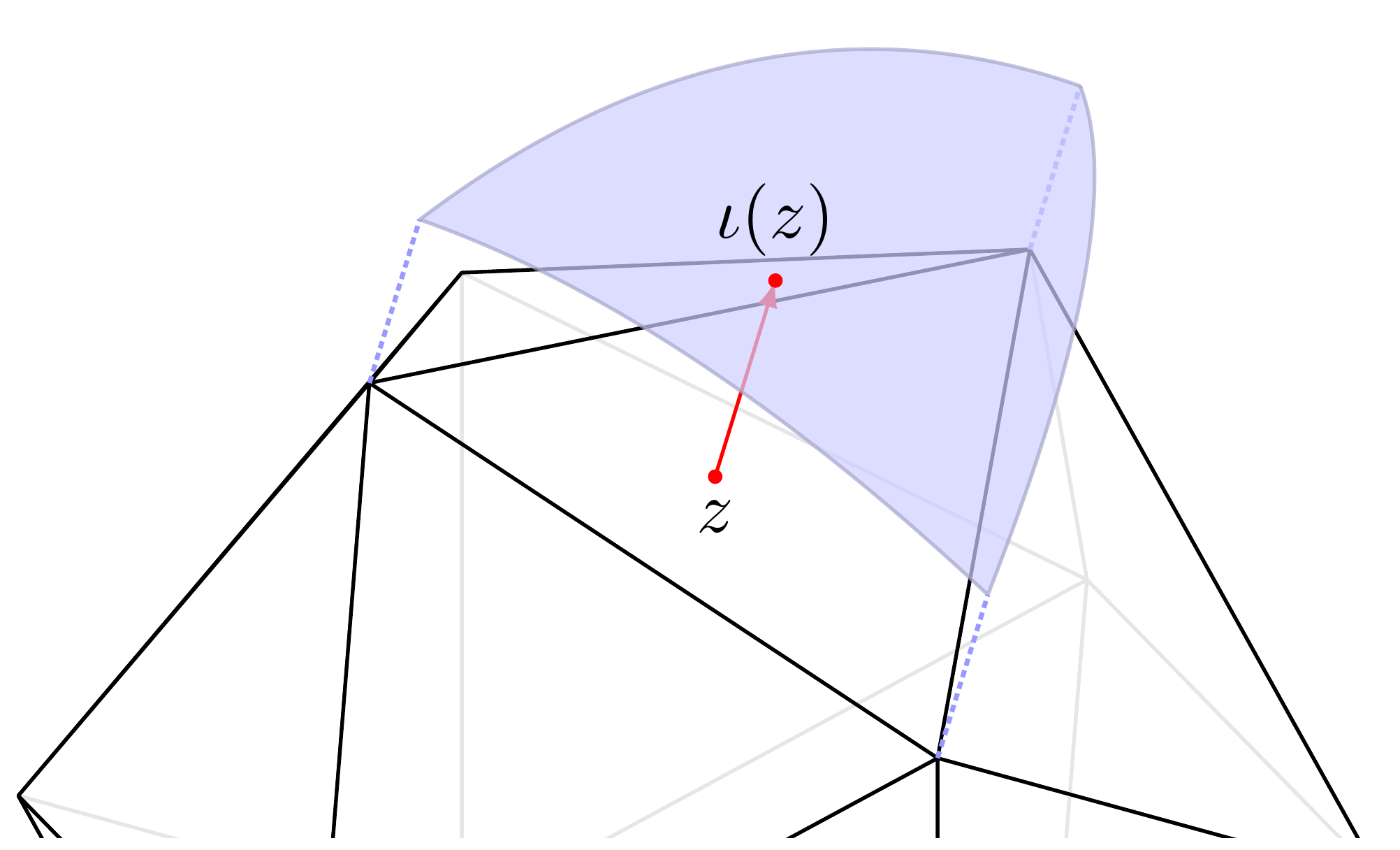}
    \caption{
        Each simplex $T$ in a triangulation (black wireframe plot) is in
        homeomorphic correspondence to a piece $\iota(T)$ of the original
        manifold (blue) through the map $\iota\colon \IM_h \to \IM$.
    }\label{fig:trisurf-homeo_vogt}
\end{figure}

\begin{figure}[t]
\sidecaption[t]
    \includegraphics[width=0.64\textwidth]{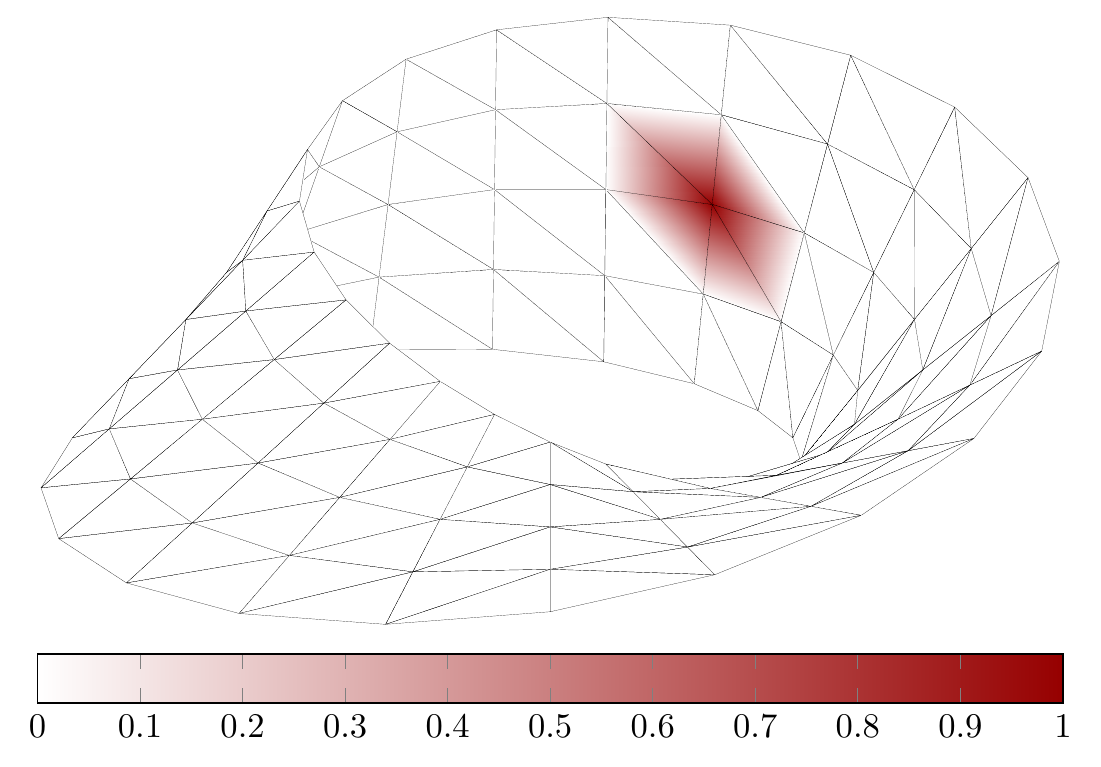}
    \caption{
        The first-order finite element space $S_h$ is spanned by a nodal basis
        $\chi_1,\dots,\chi_L$ which is uniquely determined by the property
        $\chi_k(Z^l) = 1$ if $k=l$ and $\chi_k(Z^l) = 0$ otherwise.
        The illustration shows a triangulation of the Moebius strip with a color
        plot of a nodal basis function.
    }\label{fig:moebius-basis_vogt}
\end{figure}

We translate the finite element approach for functional lifting proposed in
\cite{mollenhoff2017_vogt} to the manifold-valued setting by employing the
notation from surface finite element methods \cite{dziuk2013_vogt}.

The manifold $\IM \subset \R^N$ is approximated by a triangulated
topological manifold $\IMh \subset \R^N$ in the sense that
there is a homeomorphism $\iota\colon \IMh \to \IM$ (Fig.~\ref{fig:trisurf_vogt}
and \ref{fig:trisurf-homeo_vogt}).
By $\ITh$, we denote the set of simplices that make up $\IMh$:
\begin{equation}
    \bigcup_{T \in \ITh} T = \IMh.
\end{equation}
For $T,\tilde{T} \in \ITh$, either $T \cap \tilde{T} = \emptyset$ or
$T \cap \tilde{T}$ is an ($s-k$)-dimensional face 
for
$k \in \{1,\dots,s\}$.
Each simplex $T \in \ITh$ spans an $s$-dimensional linear subspace
of $\R^N$ and there is an orthogonal basis representation $P_T \in \R^{s,N}$ of
vectors in $\R^N$ to that subspace.
Furthermore, for later use, we enumerate the vertices of the triangulation as
$Z^1,\dots,Z^L \in \IM \cap \IM_h$.

For the numerics, we assume the first-order finite element space
\begin{equation}
    S_h := \{ \phi_h \in C^0(\IMh) :
        \phi_h|_T \text{ is linear affine for each } T \in \ITh
    \}.
\end{equation}
The functions in $S_h$ are piecewise differentiable on $\IMh$ and
we define the surface gradient $\nabla_T \phi_h \in \R^{N,d}$  of
$\phi_h \in S_h$ by the gradient of the linear affine extension of $\phi_h|_T$
to $\R^N$.
If $L$ is the number of vertices in the triangulation of $\IM_h$, then
$S_h$ is a linear space of dimension $L$ with nodal basis $\chi_1,\dots,\chi_L$
which is uniquely determined by the property $\chi_k(Z^l) = 1$ if $k=l$ and
$\chi_k(Z^l) = 0$ otherwise (Fig.~\ref{fig:moebius-basis_vogt}).

The dual space of $S_h$, which we denote by $\IfM_h(\IM_h)$, is a space of
signed measures.
We identify $\IfM_h(\IM_h) = \R^L$ via dual pairing with the nodal basis
$\chi_1,\dots,\chi_L$, i.e., to each $\mu_h \in \IfM_h(\IM_h)$ we associate the
vector $(\langle \mu_h, \chi_1 \rangle, \dots, \langle \mu_h, \chi_L \rangle)$.
We then replace the space $\IP(\IM)$ of probability measures over $\IM$
by the convex subset
\begin{equation}
    \IP_h(\IM_h) = \left\{ \mu_h \in \IfM_h(\IM_h) :
        \mu_h \geq 0,
        \sum_{k=1}^L \langle \mu_h, \chi_k \rangle = 1
    \right\}.
\end{equation}

The energy functional is then translated to the discretized setting by
redefining the integrand $f$ on $\IMh$ for any $x \in \Omega$, $z \in \IMh$
and $\xi \in \R^{N,d}$ as
\begin{equation}
    \tilde f(x,z,\xi) := \rho(x,\iota(z)) + \eta(P_T\xi)
\end{equation}
The epigraphical constraints in $\IK$ translate to
\begin{equation}\label{eq:epi-constr-discr_vogt}
    \forall x \in \Omega \forall z \in \IM_h\colon\quad
        \eta^*(P_T \nabla_z p(x,z)) - \rho(x,\iota(z)) + q(x,z) \leq 0,
\end{equation}
for functions $p \in S_h^{d}$ and $q \in S_h$.
The constraints can be efficiently implemented on each $T \in \ITh$
where $\nabla_z p$ is constant and 
$q(x,z) = \langle q_{T,1}(x), z \rangle + q_{T,2}(x)$ is linear affine in $z$:
\begin{equation}
        \eta^*(P_T \nabla_T p(x))
        + \langle q_{T,1}(x), z \rangle - \rho(x,\iota(z)) \leq -q_{T,2}(x),\label{eq:epicons-raw}
\end{equation}
for any $x \in \Omega$, $T \in \ITh$ and $z \in T$.
Following the approach in \cite{mollenhoff2017_vogt}, we define 
\begin{equation}
    \rho_T^*(x,z) := \sup_{z' \in T} \langle z, z' \rangle - \rho(x,\iota(z')),
\end{equation}
and introduce auxiliary variables $a_T,b_T$ to split the epigraphical constraint \eqref{eq:epicons-raw} into two epigraphical and one linear
constraint for $x \in \Omega$ and $T \in \ITh$:
\begin{align}
    \eta^*(P_T \nabla_T p(x)) &\leq a_T(x), \\
    \rho_T^*(q_{T,1}(x)) &\leq b_T(x), \\
    a_T(x) + b_T(x) &= -q_{T,2}(x).
\end{align}
The resulting optimization problem is described by the following saddle point
form over functions $v\colon \Omega \to \IP_h(\IM_h)$,
$p \in C^1(\Omega, S_h^{d+1})$ and $q \in C(\Omega,S_h)$:
\begin{align}
    \inf_{v} \,\sup_{p,q}~~
        &\int_\Omega \langle -\Div_x p(x,\cdot) + q(x,\cdot), v(x) \rangle\,dx \\
    \text{subject to}~~
        &\eta^*(P_T \nabla_T p(x)) \leq a_T(x), \\
        &\rho_T^*(q_{T,1}(x)) \leq b_T(x), \\
        &a_T(x) + b_T(x) + q_{T,2}(x) = 0.
\end{align}
Finally, for the fully discrete setting, the domain $\Omega$ is replaced by a
Cartesian rectangular grid with finite differences operator $\nabla_x$ and
Neumann boundary conditions.

\subsection{Relation to \cite{lellmann2013_vogt}}

\begin{figure}[t]
\includegraphics[width=0.49\textwidth,trim={0 0 330 0},clip]{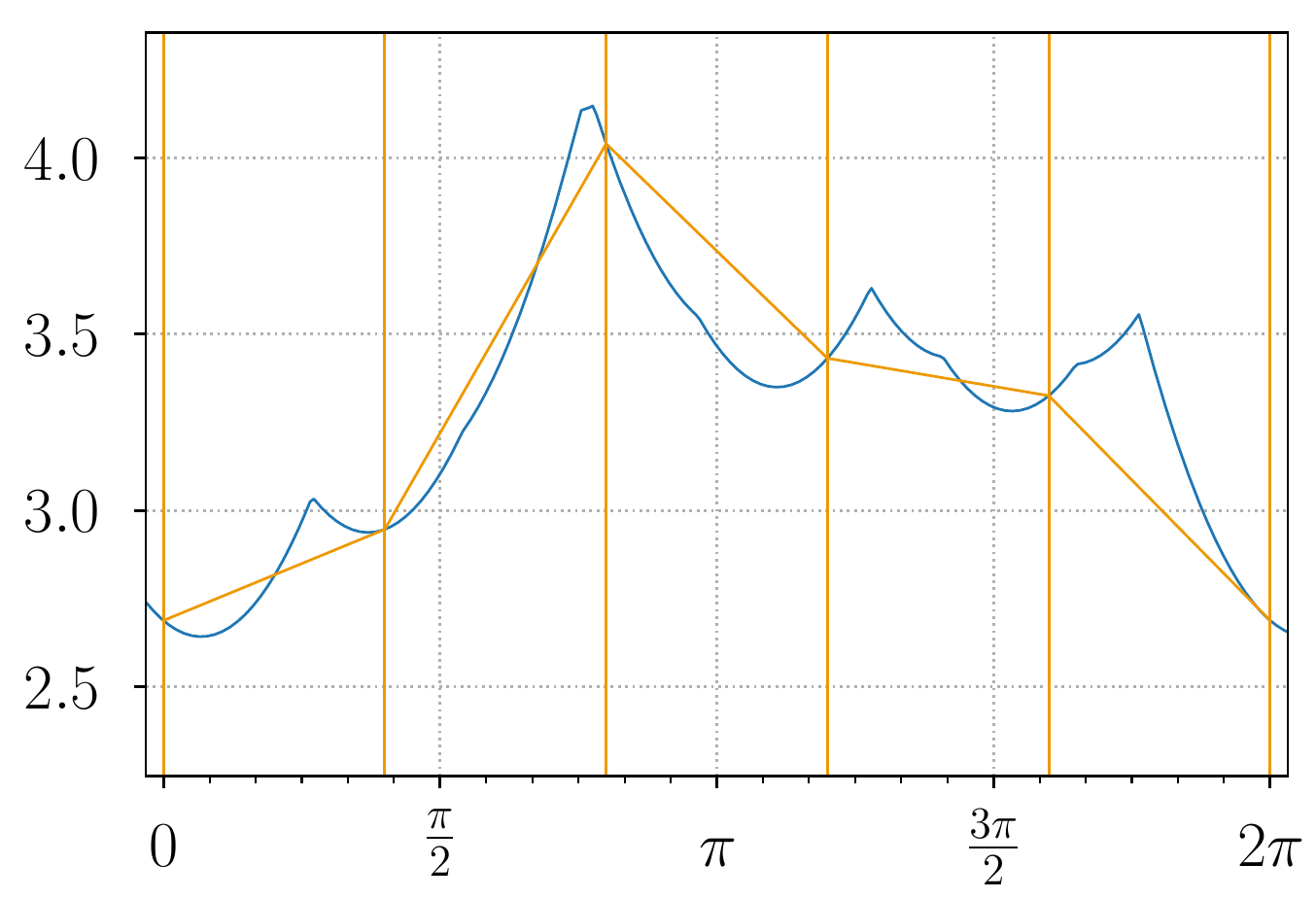}
\hfill
\includegraphics[width=0.49\textwidth,trim={0 0 330 0},clip]{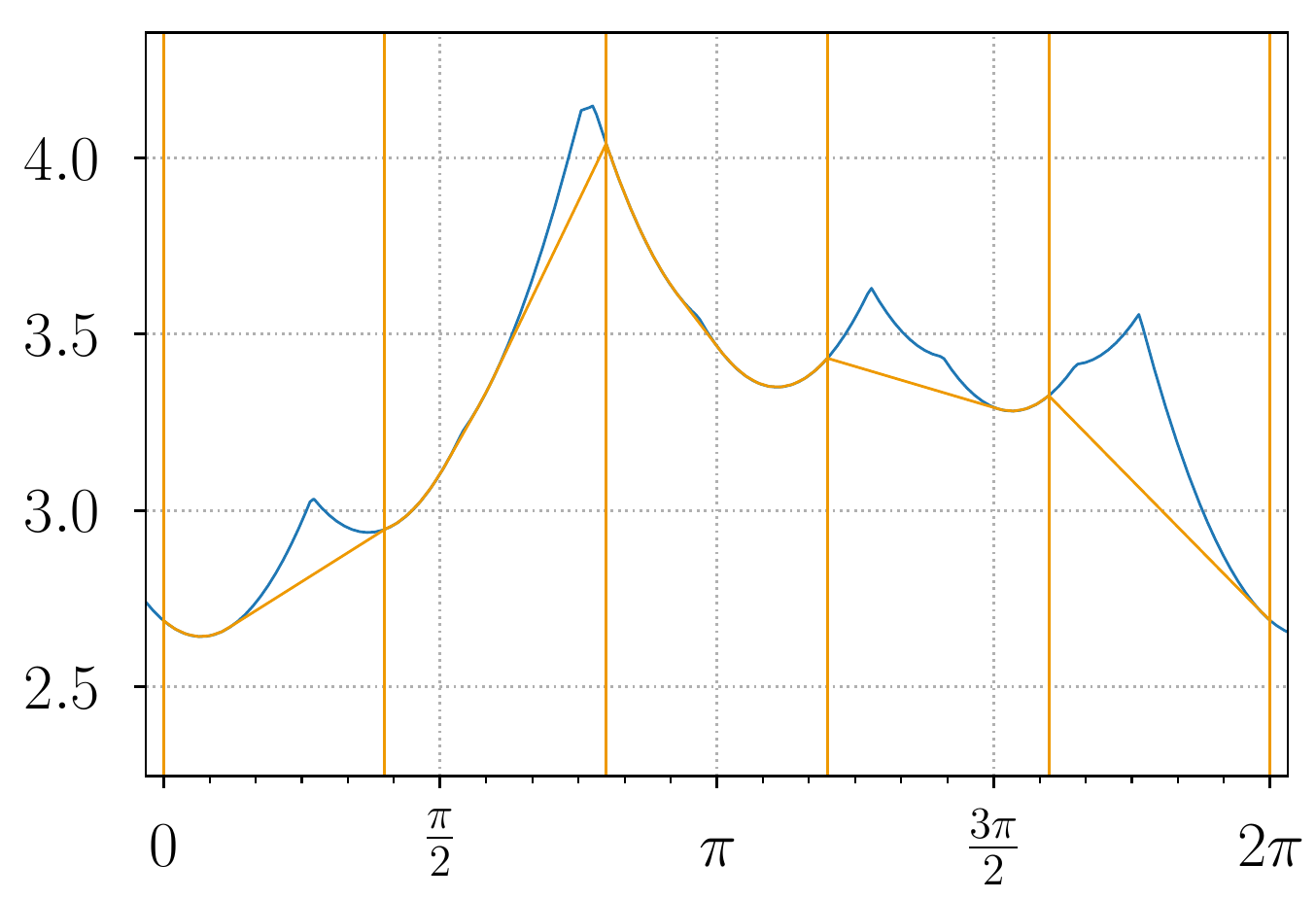}
\caption{%
    Data term discretization for the lifting approach applied to the Riemannian
    center of mass problem introduced in Figure~\ref{fig:circle-riem-com_vogt}.
    For each $x \in \Omega$, the data term $z \mapsto \rho(x,z)$ (blue graph)
    is approximated (orange graphs) between the label points $Z^k$
    (orange vertical lines).
    \textbf{Left:} In the lifting approach \cite{lellmann2013_vogt} for
    manifold-valued problems, the data term is interpolated linearly
    between the labels.
    \textbf{Right:} Based on ideas from recent scalar and vectorial lifting
    approaches \cite{mollenhoff2016_vogt,laude2016_vogt}, we interpolate
    piecewise convex between the labels.
}\label{fig:circle-rcom-sublabels_vogt}
\end{figure}

\begin{figure}[t]
\centering
\begin{tikzpicture}[
    scale=0.4,
    circle node/.style={circle,fill=black,inner sep=0,minimum size=#1}
]
    \setlength{\crvogt}{5cm}
    \node at (-30:{1.15\crvogt}) {$\IS^1$};
    \draw (0,0) circle (\crvogt);
    \draw[thick,red!90!black] ( 0:\crvogt) -- node[below] {$T$} (110:\crvogt);
    \draw[thick,red!90!black] (55:\crvogt)
        -- ++(145:{-pi*55/180*\crvogt}) node[circle node=2pt,label=30:{$v_T^1$}] {};
    \draw[thick,red!90!black] (55:\crvogt)
        -- ++(145:{pi*55/180*\crvogt}) node[circle node=2pt,label=50:{$v_T^2$}] {};
    \draw (110:\crvogt) -- (170:\crvogt);
    \draw (170:\crvogt) -- (230:\crvogt);
    \draw (230:\crvogt) -- (300:\crvogt);
    \draw (300:\crvogt) -- (360:\crvogt);
    \node[circle node=3pt,label= 55:{$y_T$}] at ( 55:\crvogt) {};
    \node[circle node=3pt,label=-20:{$Z_T^1$}] at (  0:\crvogt) {};
    \node[circle node=3pt,label=100:{$Z_T^2$}] at (110:\crvogt) {};
\end{tikzpicture}
\caption{%
    Mapping a simplex $T$, spanned by $Z_T^1,\dots,Z_T^{s+1}$, to the
    tangent space at its center-of-mass $y_T$ using the logarithmic map.
    The proportions of the simplex spanned by the mapped points
    $v_T^1,\dots,v_T^{s+1}$ may differ from the proportions of the original
    simplex for curved manifolds.
    The illustration shows the case of a circle $\IS^1 \subset \R^2$, where
    the deformation reduces to a multiplication by a scalar $\alpha_T$, the
    ratio between the geodesic (angular) and Euclidean distance between $Z_T^1$
    and $Z_T^2$.
    The gradient $\nabla_T p$ of a finite element $p \in S_h$ can be modified
    according to this change in proportion in order to make up for some of the
    geometric (curvature) information lost in the discretization.
}\label{fig:circle-gradient_vogt}
\end{figure}
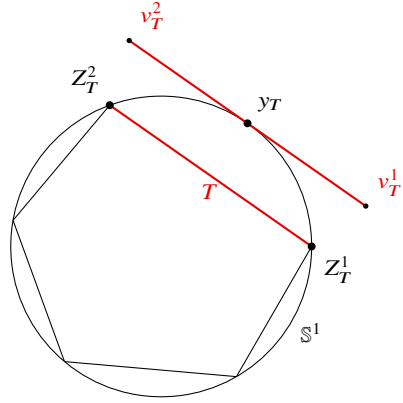

In \cite{lellmann2013_vogt}, a similar functional lifting is proposed for
the special case of total variation regularization and without the finite
elements interpretation. More precisely, the regularizing term is chosen to be
$\eta(\xi) = \lambda \|\xi\|_{\sigma,1}$ for $\xi \in \R^{s,d}$, where
$\|\cdot\|_{\sigma,1}$ is the matrix nuclear norm, also
known as Schatten-$1$-norm, which is given by the sum of singular values of a
matrix.
It is the dual to the matrix operator or spectral norm
$\|\cdot\|_{\sigma,\infty}$.
If we substitute this choice of $\eta$ into the discretization given above, the
epigraphical constraint \eqref{eq:epi-constr-discr_vogt} translates to the two
constraints
\begin{equation}
    \|P_T \nabla_T p(x)\|_{\sigma,\infty} \leq \lambda
        \text{ and }
    q(x,z) \leq \rho(x,\iota(z)).
\end{equation}
The first one is a Lipschitz constraint just as in the model from
\cite{lellmann2013_vogt}, but two differences remain:

\begin{enumerate}
\item In \cite{lellmann2013_vogt}, the lifted and discretized form of the data
    term reads
    \begin{equation}
        \int_\Omega \sum_{k=1}^L \rho(x,Z^k) v(x)^k \,dx.
    \end{equation}
    This agrees with our setting if $z \mapsto \rho(x,\iota(z))$ is
    affine linear on each simplex $T \in \ITh$, as then
    $q(x,z) = \rho(x,\iota(z))$ maximizes the objective function for any $p$
    and $v$.
    Hence, the model in \cite{lellmann2013_vogt} doesn't take into account 
    any information about $\rho$ below the resolution of the triangulation.
    We improve this by implementing the epigraph constraints
    $\rho_T^*(q_{T,1}(x)) \leq b_T(x)$ as proposed in \cite{laude2016_vogt}
    using a convex approximation of $\rho_T$ (see
    Fig.~\ref{fig:circle-rcom-sublabels_vogt}).
    The approximation is implemented numerically with piecewise affine linear
    functions in a ``sublabel-accurate'' way, i.e., at a resolution below 
    the resolution of the triangulation .
\item A very specific discretization of the gradients $\nabla_T p(x)$ is
    proposed in \cite{lellmann2013_vogt}:
    To each simplex in the triangulation a mid-point $y_T \in \IM$ is
    associated.
    The vertices $Z_T^1, ...,Z_T^{s+1}$ of the simplex are projected to the
    tangent space at $y_T$ as $v_T^k := \log_{y_T} Z_T^k$.
    The gradient is then computed as the vector $g$ in the tangent space
    $T_{y_T}\IM$ describing the affine linear map on $T_{y_T}\IM$ that takes
    values $p(Z_T^k)$ at the points $v_T^k$, $k=1,\dots,s+1$.
    
    This procedure aims to make up for the error introduced by the simplicial
    discretization and amounts to a different choice of $P_T$ -- a slight variant of our model.
    We did not observe any significant positive or negative effects
    from using either discretization; the
    difference between the minimizers is very small.
    
    In the one-dimensional case, the two approaches differ only in a constant factor:
    Denote by $P_T \in \R^{s,N}$ the orthogonal basis representation of
    vectors in $\R^N$ in the subspace spanned by the simplex $T \in \ITh$
    and denote by $\tilde P_T \in \R^{s,N}$ the alternative approach from
    \cite{lellmann2013_vogt}.
    Now, consider a triangulation $\ITh$ of the circle $\IS^1 \subset \R^2$ and
    a one-dimensional simplex $T \in \ITh$.
    A finite element $p \in S_h$ that takes values $p_1,p_2 \in \R$ at the vertices
    $Z_T^1, Z_T^2 \in \R^2$ that span $T$ has the gradient
    \begin{equation}
        \nabla_T p = (p_1 - p_2) \frac{Z_T^1 - Z_T^2}{\|Z_T^1 - Z_T^2\|_2^2}
            \in \R^{2}
    \end{equation}
    and $P_T, \tilde P_T \in \R^{1,2}$ are given by
    \begin{align}
        P_T &:= \frac{(Z_T^1 - Z_T^2)^\top}{\|Z_T^1 - Z_T^2\|_2},\qquad
        \tilde P_T := \frac{(Z_T^1 - Z_T^2)^\top}{d_{\IS^1}(Z_T^1,Z_T^2)}.
    \end{align}
    Hence $P_T = \alpha_T \tilde P_T$ for
    $\alpha_T = d_{\IS^1}(Z_T^1,Z_T^2)/\|Z_T^1 - Z_T^2\|_2$
    the ratio between geodesic (angular) and Euclidean distance between the
    vertices.
    If the vertices are equally spaced on $\IS^1$, this is a constant factor
    independent of $T$ that typically scales the discretized regularizer by
    a small constant factor.
    On higher-dimensional manifolds, more general linear transformations
    $P_T = A_T \tilde P_T$ come into play.
    For very irregular triangulations and coarse discretization, this may affect the minimizer; however, in our experiments the observed differences were negligible.
\end{enumerate}

\subsection{Full discretization and numerical implementation}
A prime advantage of the lifting method when applied to manifold-valued problems is that it translates most parts of the problem into Euclidean space. This allows to apply established solution strategies for the non-manifold case, which rely on non-smooth convex optimization:
After discretization, the convex-concave saddle-point form allows for a solution
using the primal-dual hybrid gradient method
\cite{chambolle2011_vogt,chambolle2012_vogt} with recent extensions
\cite{goldstein2013_vogt}.
In this optimization framework, the epigraph constraints are realized by
projections onto the epigraphs in each iteration step.
For the regularizers to be discussed in this paper (TV, quadratic and Huber), we refer
to the instructions given in \cite{pock2010_vogt}.
For the data term $\rho$, we follow the approach in
\cite{laude2016_vogt}:
For each $x \in \Omega$, The data term $z \mapsto \rho(x,\iota(z))$ is sampled
on a subgrid of $\IM_h$ and approximated by a piecewise affine linear function.
The quickhull algorithm can then be used to get the convex hull of this
approximation.
Projections onto the epigraph of $\rho^*_T$ are then projections
onto convex polyhedra, which amounts to solving many low-dimensional
quadratic programs; see \cite{laude2016_vogt} for more details.

Following \cite{lellmann2013_vogt}, the numerical solution
$u\colon \Omega \to \IP_h(\IM_h)$, taking values in the lifted space
$\IP_h(\IM_h)$, is projected back to a function $u\colon \Omega \to \IM$, taking
values in the original space $\IM$,
by mapping, for each $x \in \Omega$ separately, a probability measure
$u(x) = (\lambda_1,\dots,\lambda_L) = \mu_h \in \IP_h(\IM_h)$
to the following Riemannian center of mass on the original manifold $\IM$:
\begin{equation}
    \mu_h = (\lambda_1,\dots,\lambda_L)
        \mapsto
    \argmin_{z \in \IM} \sum_{k=1}^L \lambda_k d_\IM(z,Z^k)^2
\end{equation}
For $\IM = \R^s$, this coincides with the usual weighted mean
$\bar{z} = \sum_{k=1}^L \lambda_k Z_k$.
However, on manifolds this minimization is known to be a non-convex problem
with non-unique solutions (compare Fig.~\ref{fig:circle-riem-com_vogt}).
Still, in practice the iterative method described in \cite{karcher1977_vogt}
yields reasonable results for all real-world data considered in this work:
Starting from a point $z_0 := Z^k$ with maximum weight $\lambda_k$, we proceed
for $i \geq 0$ by projecting the $Z^k$, $k=1,\dots,L$, to the tangent space at
$z_i$ using the inverse exponential map, taking the linear weighted mean $v_i$
there and defining $z_{i+1}$ as the projection of $v_i$ to $\IM$ via the
exponential map:
\begin{align}
    V_i^k &:= \log_{z_i}(Z^k) \in T_{z_i}\IM,~ k=1,\dots,L, \\
    v_i &:= \sum_{k=1}^L \lambda_k V_i^k \in T_{z_i}\IM, \\
    z_{i+1} &:= \exp_{z_i}(v_i).
\end{align}
The method converges rapidly in practice. It has to be applied only once for each $x \in \Omega$ after solving the lifted problem, so that
efficiency is non-critical.


%% file: chpt/sec2.tex
\section{Numerical Results}

We apply our model to problems with quadratic data term
$\rho(x,z) := d_\IM^2(I(x),z)$ and Huber, total variation (TV) and Tikhonov
(quadratic) regularization with parameter $\lambda > 0$:
\begin{align}
    \eta_{\text{TV}}(\xi) &:= \lambda \|\xi\|_2, \\
    \eta_{\text{Huber}}(\xi) &:= \lambda \phi_\alpha(\xi), \\
    \eta_{\text{quad}}(\xi) &:= \frac{\lambda}{2} \|\xi\|_2^2,
\end{align}
where the Huber function $\phi_\alpha$ for $\alpha > 0$ is defined by
\begin{equation}
    \phi_\alpha(\xi) := \begin{cases}
        \frac{\|\xi\|_2^2}{2\alpha} & \text{if $\|\xi\|_2 \leq \alpha$,} \\
        \|\xi\|_2 - \frac{\alpha}{2} & \text{if $\|\xi\|_2 > \alpha$.}
    \end{cases}
\end{equation}
Note that previous lifting approaches for manifold-valued data were restricted to total variation regularization $\eta_{\text{TV}}$.

The methods were implemented in Python 3 with NumPy and PyCUDA, running on an Intel Core i7 4.00\,GHz with an NVIDIA GeForce GTX 1080 Ti 12 GB and 16\,GB RAM. The iteration was stopped as soon as the relative gap between primal and dual objective fell below $10^{-5}$. Approximate runtimes ranged between 5 and 45 minutes.
The code is available from \url{https://github.com/room-10/mfd-lifting}.

\subsection{One-dimensional denoising on a Klein bottle}

\begin{figure}[t]
    \includegraphics[angle=-90,width=\textwidth,trim={0 0 0 0},clip]{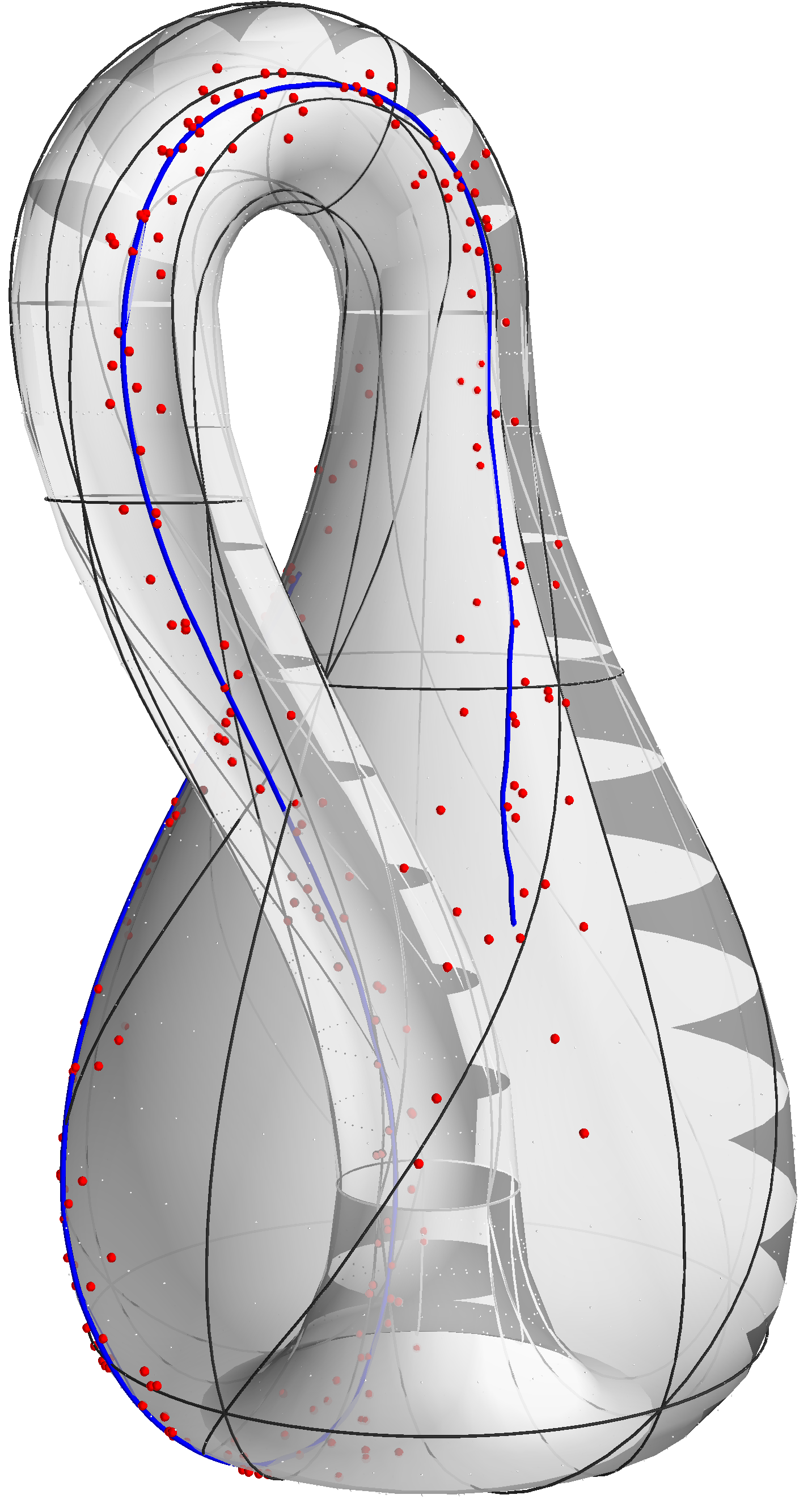}
    \caption{
        Tikhonov (quadratic) denoising (blue) of a
        one-dimensional signal (red) $u\colon [0,1] \to \IM$ with values on
        the two-dimensional Klein surface (commonly referred to as Klein bottle)
        $\IM \subset \R^3$.
        The black wireframe lines on the surface represent the triangulation
        used by the discretization of our functional lifting
        approach.
        The numerical implementation recovers the denoised signal at a
        resolution far below the resolution of the manifold's discretization.
        The lifting approach does not require the
        manifold to be orientable.
    }\label{fig:klein_vogt}
\end{figure}

Our model can be applied to both orientable and non-orientable manifolds.
Figure~\ref{fig:klein_vogt} shows an application of our method to Tikhonov
denoising of a synthetic one-dimensional signal $u\colon [0,1] \to \IM$ on the
two-dimensional Klein surface embedded in $\R^3$, a non-orientable closed
surface that cannot be embedded into $\R^3$ without self-intersections.
Our numerical implementation uses a triangulation with a very low count of $5 \times 5$ vertices
and 50 triangles.
The resolution of the signal (250 one-dimensional data points) is far below the
resolution of the triangulation and, still, our approach is able to restore
a smooth curve.

\subsection{Three-dimensional manifolds: $SO(3)$}

\begin{figure}[t]
    \includegraphics[width=\textwidth,trim={150  60 110 110},clip]{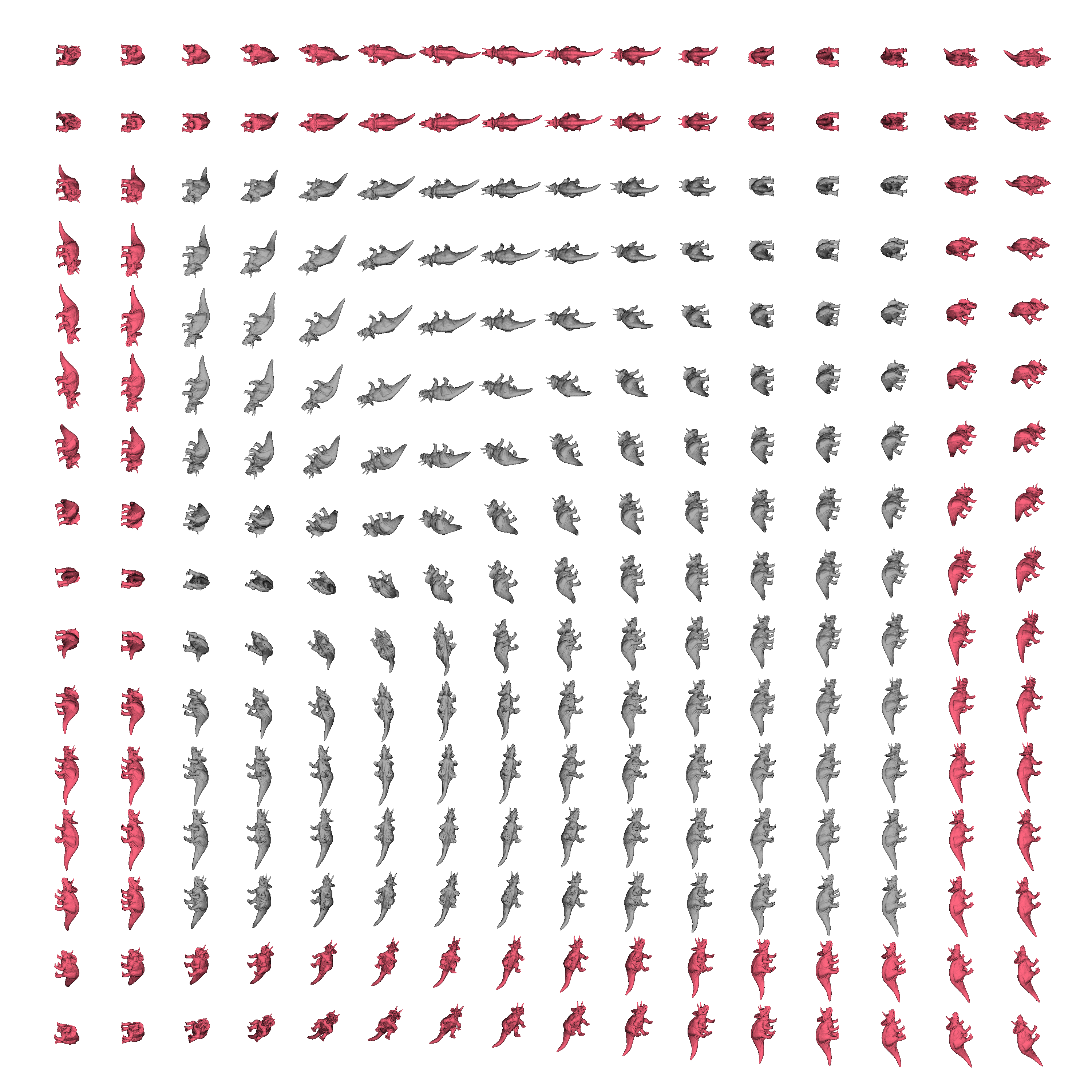}
    \caption{
        Tikhonov inpainting of a 2-dimensional signal of
        (e.g., camera) orientations, elements of the three-dimensional special
        orthogonal group of rotations $SO(3)$, a manifold of dimension
        $s=3$.
        The masked input signal (red) is inpainted (gray) using our model
        with Tikhonov (quadratic) regularization.
        The interpolation into the central area is smooth.
        Shape: \emph{Triceratops} by BillyOceansBlues (CC-BY-NC-SA,
        \url{https://www.thingiverse.com/thing:3313805}).
    }\label{fig:camera_vogt}
\end{figure}

Signals with rotational range $u\colon \Omega \to SO(3)$ occur in the
description of crystal symmetries in EBSD (Electron Backscatter Diffraction Data)
and in motion tracking. 
The rotation group $SO(3)$ is a three-dimensional manifold that can be identified
with the three-dimensional unit-sphere $\IS^3$ up to identification of
antipodal points via the quaternion representation of 3D rotations.
A triangulation of $\IS^3$ is given by the vertices and simplicial faces of the
hexacosichoron (600-cell), a regular polytope in $\R^4$ akin to the icosahedron
in $\R^3$.
As proposed in \cite{lellmann2013_vogt}, we eliminate opposite points in the
hexacosichoron and obtain a discretization of $SO(3)$ with $60$ vertices and
$300$ tetrahedral faces.

Motivated by B\'{e}zier surface interpolation \cite{absil2016_vogt}, we applied
Tikhonov regularization to a synthetic inpainting (interpolation) problem with added noise (Fig.~\ref{fig:camera_vogt}).
In our variational formulation, we chose $\rho(x,z) = 0$ for $x$ in the
inpainting area and $\rho(x,z) = \delta_{\{z=I(x)\}}$ (a hard constraint to the
input signal $I\colon \Omega \to SO(3)$) for $x$ in the known area.

Using the proposed sublabel-accurate handling of data terms, we obtain good
results with only $60$ vertices, in contrast to \cite{lellmann2013_vogt},
where the discretization is refined to $720$ vertices (Fig.~\ref{fig:camera_vogt}).

\subsection{Normals fields from digital elevation data}

\begin{figure}
    \centering
    \includegraphics[width=0.35\textwidth,trim= 40 20 940 20,clip]{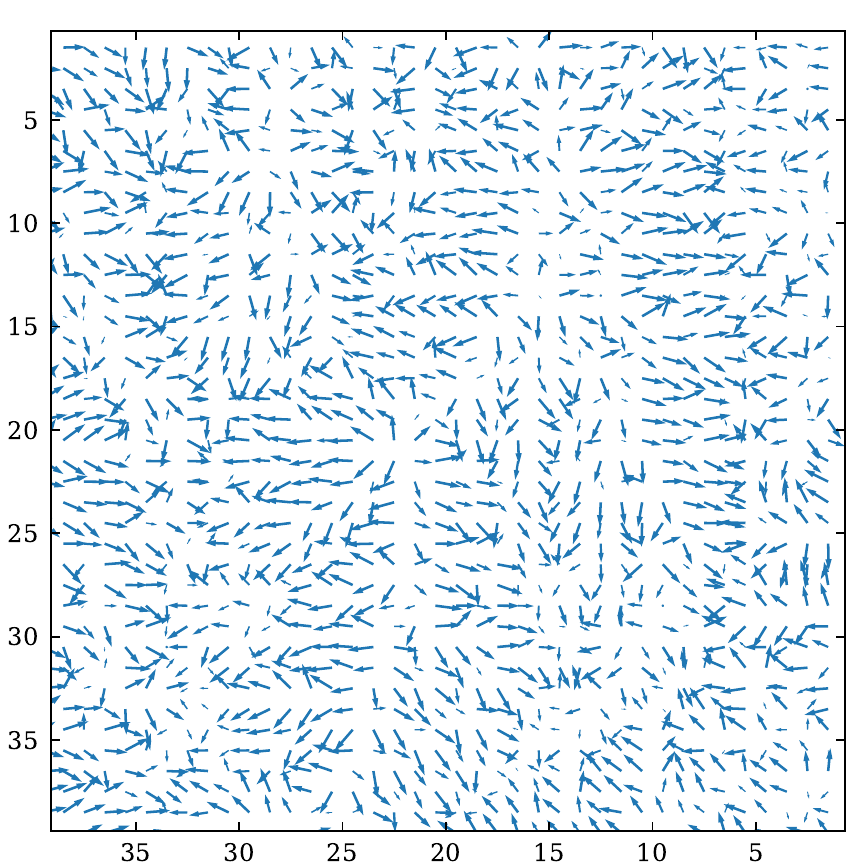}
    \hspace{0.5cm}
    \includegraphics[width=0.42\textwidth,trim=670 40 315 40,clip]{fig/bull-huber-01-075.pdf}

    \includegraphics[width=0.35\textwidth,trim=345 20 635 20,clip]{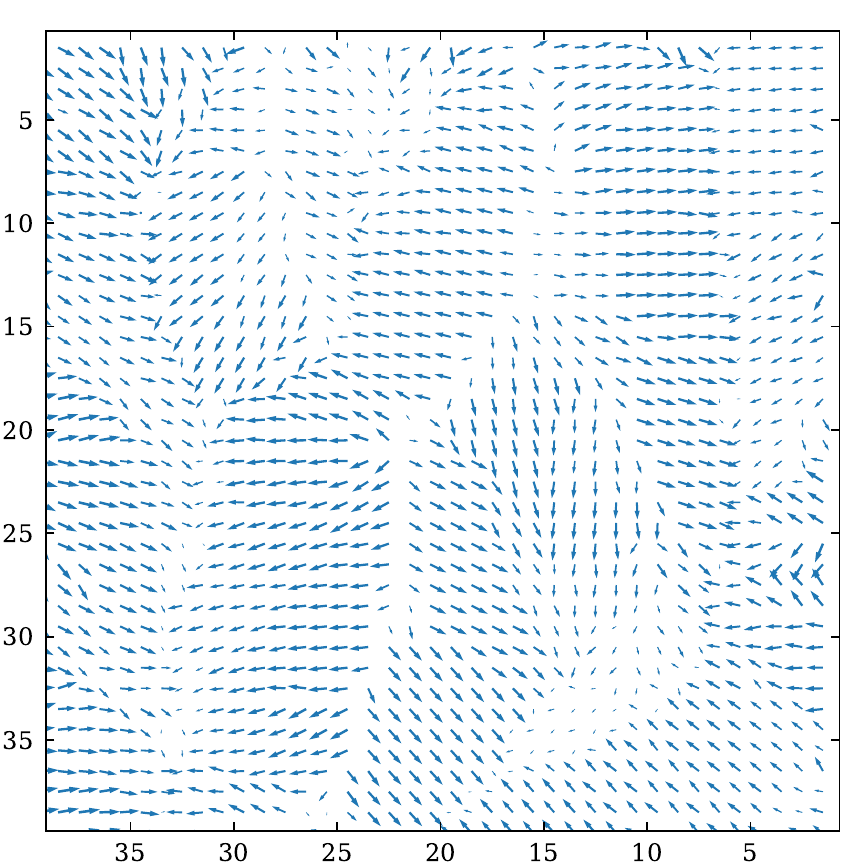}
    \hspace{0.5cm}
    \includegraphics[width=0.42\textwidth,trim=975 40  10 40,clip]{fig/bull-tv-04.pdf}

    \includegraphics[width=0.35\textwidth,trim=345 20 635 20,clip]{fig/bull-huber-01-075.pdf}
    \hspace{0.5cm}
    \includegraphics[width=0.42\textwidth,trim=975 40  10 40,clip]{fig/bull-huber-01-075.pdf}

    \includegraphics[width=0.35\textwidth,trim=345 20 635 20,clip]{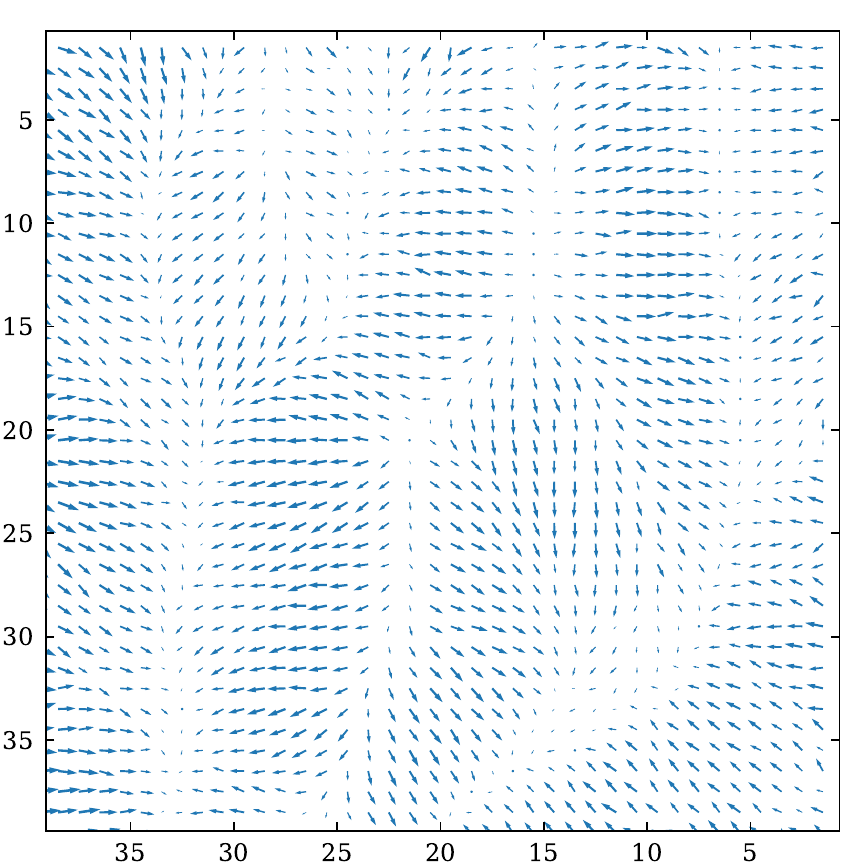}
    \hspace{0.5cm}
    \includegraphics[width=0.42\textwidth,trim=975 40  10 40,clip]{fig/bull-quadratic-3.pdf}

    \caption{
        Denoising of $\IS^2$-valued surface normals on the digital elevation model (DEM) dataset from \cite{gesch2009_vogt}:
        Noisy input (top),
        total variation ($\lambda = 0.4$) denoised image (second from top),
        Huber ($\alpha = 0.1$, $\lambda = 0.75$) denoised image (second from bottom),
        quadratically ($\lambda = 3.0$) denoised image (bottom).
        Mountain ridges are sharp while hillsides remain smooth with Huber.
        TV enforces flat hillsides and Tikhonov regularization smoothes out
        all contours.
    }\label{fig:dem_vogt}
\end{figure}

In digital elevation models (DEM), elevation information for earth science
studies and mapping applications often includes surface normals which can be
used to produce a shaded coloring of elevation maps.
Normal fields $u\colon \Omega \to \IS^2$ are defined on a rectangular image
domain $\Omega \subset \R^2$; variational processing of the normal fields is therefore a manifold-valued problem on the two-dimensional sphere
$\IS^2 \subset \R^3$.

Denoising using variational regularizers from manifold-valued image processing before computing the shading
considerably improves visual quality (Fig.~\ref{fig:dem_vogt}).
For our framework, the sphere was discretized using 12 vertices and 20 triangles,
chosen to form a regular icosahedron.
The same dataset was used in \cite{lellmann2013_vogt}, where the proposed lifting
approach required 162 vertices -- and solving a proportionally larger optimization problem -- in order to produce comparable results.

We applied our approach with TV, Huber and Tikhonov regularization. Interestingly, many of the qualitative properties known from RGB and grayscale image
processing appear to transfer to the manifold-valued case:
TV enforces piecewise constant areas (flat hillsides), but preserves
edges (mountain ridges).
Tikhonov regularization gives overall very smooth results, but tends to lose
edge information.
With Huber regularization, edges (Mountain ridges) remain sharp while hillsides
are smooth, and flattening is avoided (Fig.~\ref{fig:dem_vogt}).

\subsection{Denoising of high resolution InSAR data}

\begin{figure}[t]
    \centering
    \newlength\picwidth
    \setlength\picwidth{0.496\textwidth}
    \newlength\spywidth
    \setlength\spywidth{0.25\picwidth}
    \begin{tikzpicture}[spy using outlines={%
            height = \spywidth, width = \spywidth,
            magnification = 1.7, white, connect spies
    }]
        \node[inner sep=0] (insar-noisy) at (0,0) {
            \includegraphics[
                trim= 10 0 390 0,clip,width=\picwidth
            ]{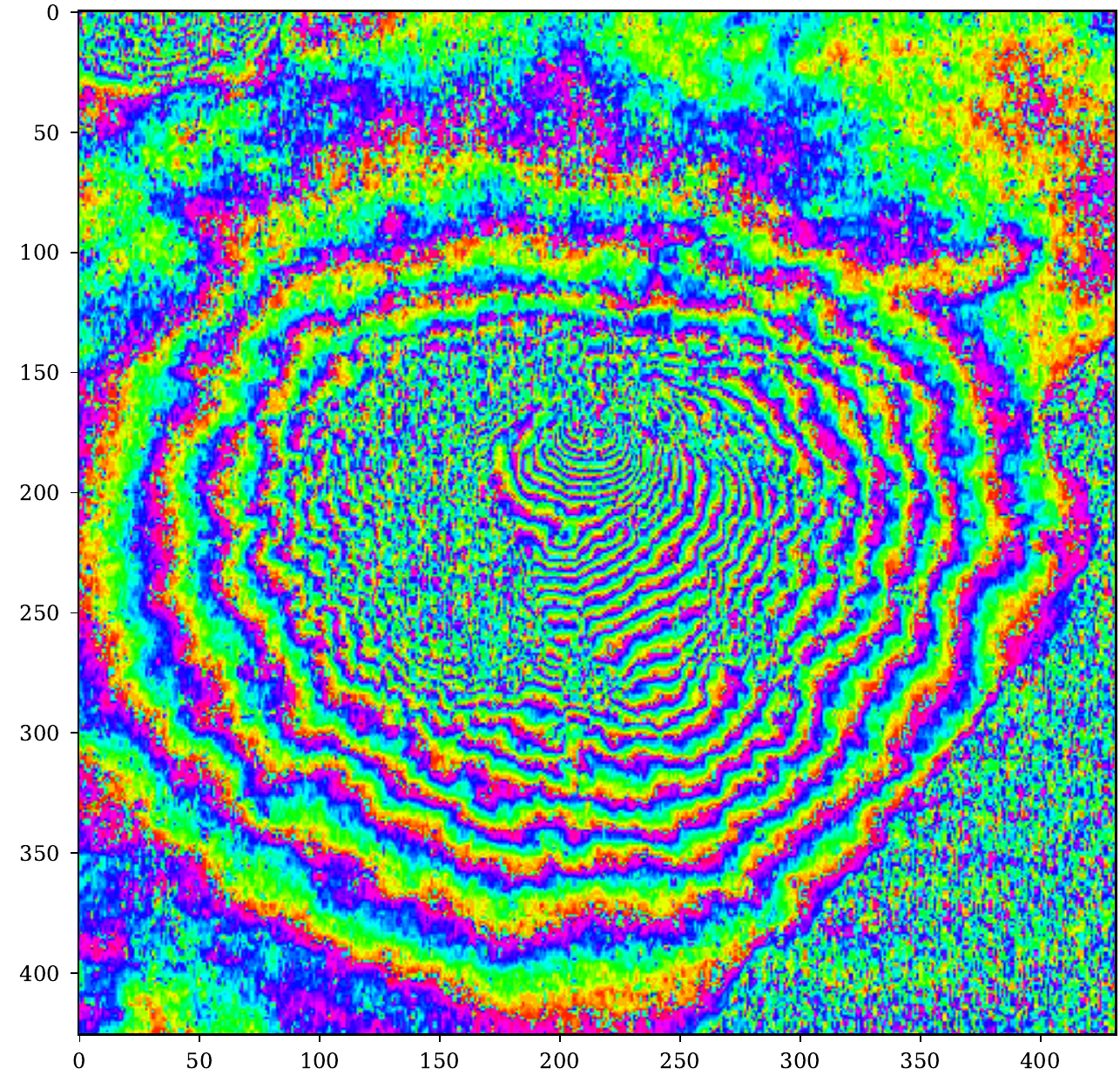}
        };
        \RelativeSpy{insar-spy1}{insar-noisy}{(.69,.86)}{(0.9,0.13)}

        \node[inner sep=0,right=0.1cm of insar-noisy] (insar-tv) {%
            \includegraphics[
                trim=390 0  10 0,clip,width=\picwidth
            ]{fig/insar-tv-06.pdf}
        };
        \RelativeSpy{insar-spy2}{insar-tv}{(.69,.86)}{(0.15,0.13)}

        \node[inner sep=0,below=0.1cm of insar-noisy] (insar-huber) {%
            \includegraphics[
                trim=390 0  10 0,clip,width=\picwidth
            ]{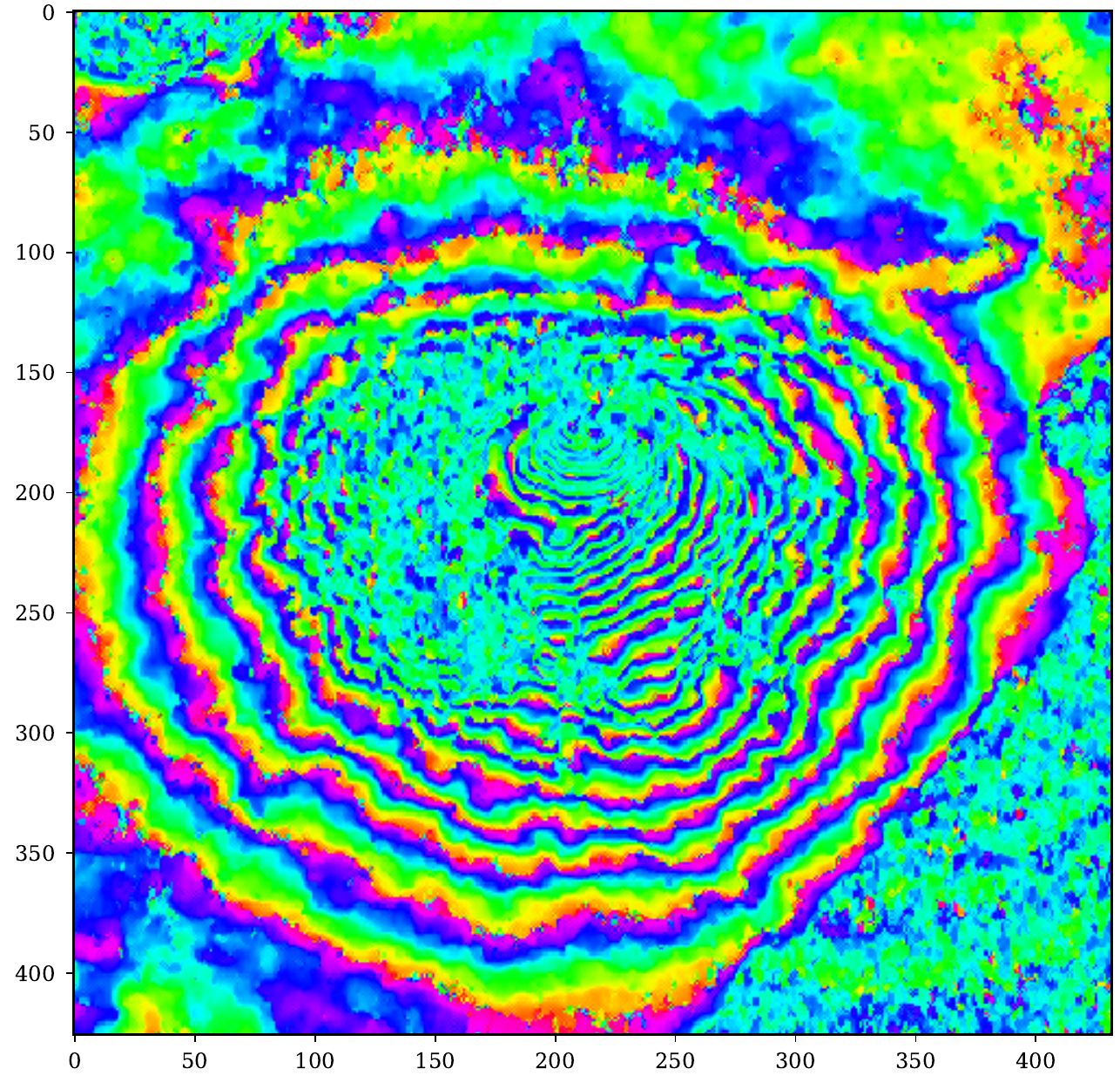}
        };
        \RelativeSpy{insar-spy3}{insar-huber}{(.69,.86)}{(0.9,0.88)}

        \node[inner sep=0,right=0.1cm of insar-huber] (insar-quadratic) {%
            \includegraphics[
                trim=390 0  10 0,clip,width=\picwidth
            ]{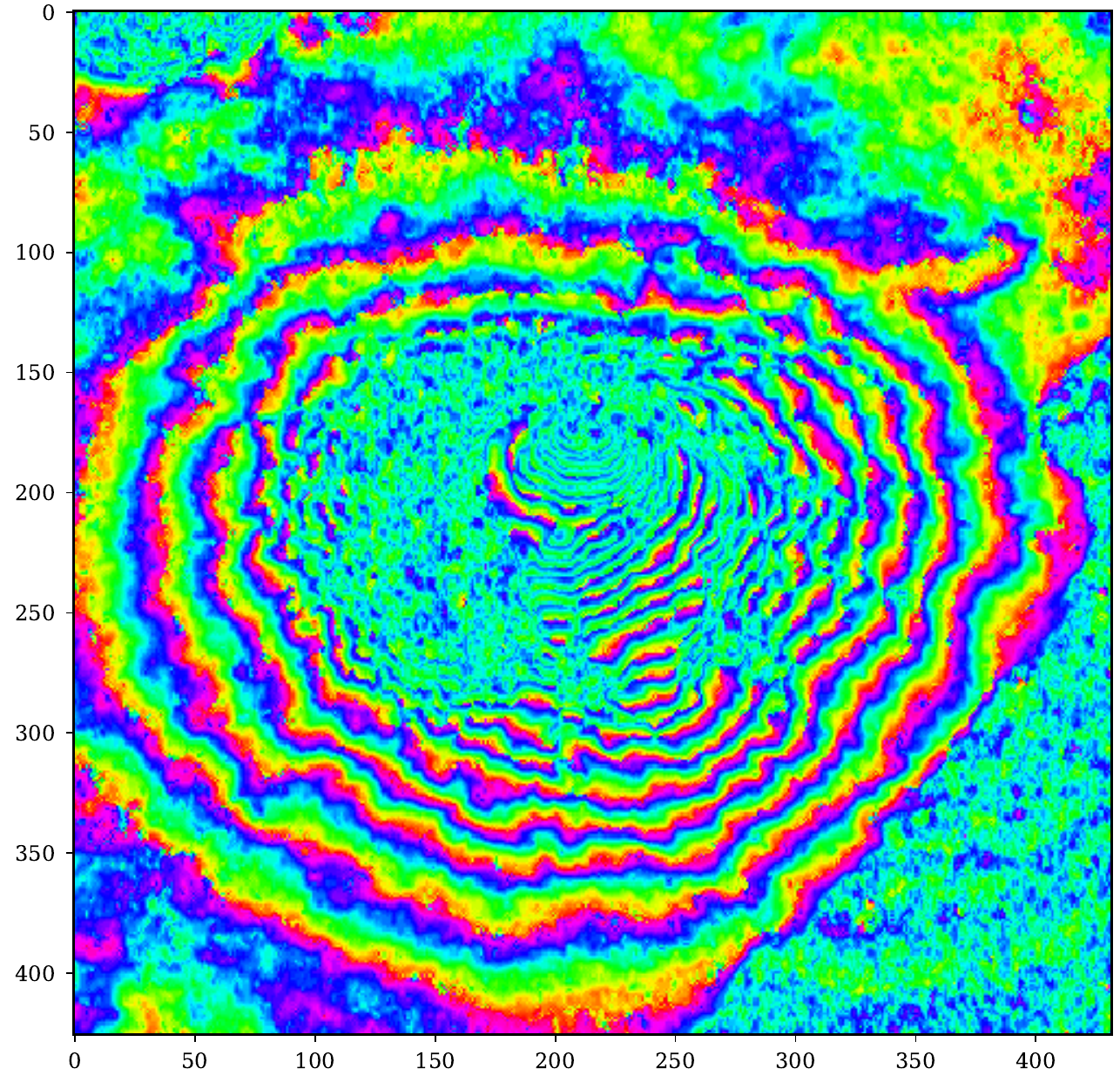}
        };
        \RelativeSpy{insar-spy4}{insar-quadratic}{(.69,.86)}{(0.15,0.88)}
    \end{tikzpicture}

    \caption{\label{fig:insar}
        Denoising of $\IS^1$-valued InSAR measurements from Mt. Vesuvius, dataset from \cite{rocca1997_vogt}:
        Noisy input (top left),
        total variation ($\lambda = 0.6$) denoised image (top right),
        Huber ($\alpha = 0.1$, $\lambda = 0.75$) denoised image (bottom left),
        quadratically ($\lambda = 1.0$) denoised image (bottom right).
        All regularization strategies successfully remove most of the noise.
        The total variation regularizer enforces clear contours, but exhibits
        staircasing effects.
        The staircasing is removed with Huber while contours are still quite
        distinct.
        Quadratic smoothing preserves some of the finer structures, but produces
        an overall more blurry and less contoured result.
    }
\end{figure}

While the resolution of the DEM dataset is quite limited ($40 \times 40$ data
points), an application to high resolution ($432 \times 426$ data points)
Interferometric Synthetic Aperture Radar (InSAR) denoising shows that our model is
also applicable in a more demanding scenario (Fig.~\ref{fig:insar}).

In InSAR imaging, information about terrain is obtained from satellite or aircraft by measuring the phase difference between the outgoing signal and the incoming reflected signal. This allows a very high relative precision, but no immediate absolute measurements, as all distances are only recovered modulo the wavelength. After normalization to $[0,2\pi)$, the phase data is correctly viewed as lying on the one-dimensional unit sphere $\IS^1$. Therefore, handling the data before any phase unwrapping is performed requires a manifold-valued framework.

Again, denoising with TV, Huber, and Tikhonov regularizations demonstrates properties comparable to those known from scalar-valued image processing while all regularization approaches reduce noise substantially (Fig.~\ref{fig:insar}).
%
%
%
%

%% file: chpt/sec3.tex
\section{Conclusion and Outlook}

We provided an overview and framework for functional lifting techniques for the variational regularization
of functions with values in arbitrary Riemannian manifolds.
The framework is motivated from the theory of currents and continuous
multi-label relaxations, but generalizes these from the context of scalar and
vectorial ranges to geometrically more challenging manifold ranges.

Using this approach, it is possible to solve variational problems for
manifold-valued images that consist of a possibly non-convex data term and
an arbitrary, smooth or non-smooth, convex first-order regularizer, such as
Tikhonov, total variation or Huber.
A refined discretization based on manifold finite element methods achieves sublabel-accurate results, which allows to use coarser discretization of the range and reduces computational effort compared to previous lifting approaches on manifolds.


A primary limitation of functional lifting methods, which equally applies to manifold-valued models, is dimensionality: The numerical cost increases exponentially with the dimensionality of the manifold due to the required discretization of the range. Addressing this issue appears possible, but will require a significantly improved discretization strategy.


%% file: vogt.bbl
\begin{thebibliography}{10}
\providecommand{\url}[1]{{#1}}
\providecommand{\urlprefix}{URL }
\expandafter\ifx\csname urlstyle\endcsname\relax
  \providecommand{\doi}[1]{DOI~\discretionary{}{}{}#1}\else
  \providecommand{\doi}{DOI~\discretionary{}{}{}\begingroup
  \urlstyle{rm}\Url}\fi

\bibitem{absil2016_vogt}
Absil, P.A., Gousenbourger, P.Y., Striewski, P., Wirth, B.: {Differentiable
  Piecewise-B{\'e}zier Surfaces on Riemannian Manifolds}.
\newblock SIAM J Imaging Sci \textbf{9}, 1788--1828 (2016)

\bibitem{absil2009_vogt}
Absil, P.A., Mahony, R., Sepulchre, R.: {Optimization algorithms on matrix
  manifolds}.
\newblock Princeton University Press (2009)

\bibitem{alberti2003_vogt}
{Alberti}, G., {Bouchitt\'e}, G., {Dal Maso}, G.: {The calibration method for
  the Mumford-Shah functional and free-discontinuity problems.}
\newblock {Calc Var Partial Differ Equ} \textbf{16}(3), 299--333 (2003)

\bibitem{bachmann2011_vogt}
Bachmann, F., Hielscher, R., Schaeben, H.: {Grain detection from 2d and 3d EBSD
  data - Specification of the MTEX algorithm}.
\newblock Ultramicroscopy \textbf{111}(12), 1720--1733 (2011)

\bibitem{bae2011_vogt}
Bae, E., Yuan, J., Tai, X.C., Boykov, Y.: {A Fast Continuous Max-Flow Approach
  to Non-convex Multi-labeling Problems}.
\newblock In: A.~Bruhn, T.~Pock, X.C. Tai (eds.) {Efficient Algorithms for
  Global Optimization Methods in Computer Vision}, pp. 134--154. Springer
  Berlin Heidelberg, Berlin, Heidelberg (2014)

\bibitem{basser1994_vogt}
Basser, P.J., Mattiello, J., LeBihan, D.: {MR diffusion tensor spectroscopy and
  imaging}.
\newblock {Biophys J} \textbf{66}(1), 259--267 (1994)

\bibitem{baust2016_vogt}
Baust, M., Weinmann, A., Wieczorek, M., Lasser, T., Storath, M., Navab, N.:
  {Combined Tensor Fitting and TV Regularization in Diffusion Tensor Imaging
  Based on a Riemannian Manifold Approach}.
\newblock IEEE Trans Med Imaging \textbf{35}, 1972--1989 (2016)

\bibitem{bacak2014_vogt}
Ba\v{c}{\'a}k, M.: {Convex Analysis and Optimization in Hadamard Spaces}.
\newblock De Gruyter (2014)

\bibitem{bacak2016_vogt}
Ba\v{c}{\'a}k, M., Bergmann, R., Steidl, G., Weinmann, A.: A {Second} {Order}
  {Nonsmooth} {Variational} {Model} for {Restoring} {Manifold}-{Valued}
  {Images}.
\newblock SIAM J Sci Comput \textbf{38}(1), A567--A597 (2016)

\bibitem{bergmann2018b_vogt}
Bergmann, R., Fitschen, J.H., Persch, J., Steidl, G.: {Priors with Coupled
  First and Second Order Differences for Manifold-Valued Image Processing}.
\newblock J Math Imaging Vis \textbf{60}, 1459--1481 (2018)

\bibitem{bergmann2018c_vogt}
Bergmann, R., Laus, F., Persch, J., Steidl, G.: {Recent Advances in Denoising
  of Manifold-Valued Images.}
\newblock Tech. Rep. arXiv:1812.08540, arXiv (2018)

\bibitem{bergmann2016_vogt}
Bergmann, R., Persch, J., Steidl, G.: {A Parallel Douglas-Rachford Algorithm
  for Minimizing ROF-like Functionals on Images with Values in Symmetric
  Hadamard Manifolds}.
\newblock SIAM J Imaging Sci \textbf{9}, 901--937 (2016)

\bibitem{bergmann2018_vogt}
Bergmann, R., Tenbrinck, D.: {A Graph Framework for Manifold-Valued Data}.
\newblock SIAM J Imaging Sci \textbf{11}, 325--360 (2018)

\bibitem{bernard2017_vogt}
Bernard, F., Schmidt, F.R., Thunberg, J., Cremers, D.: {A Combinatorial
  Solution to Non-Rigid 3D Shape-to-Image Matching}.
\newblock In: Proc ICCV 2017, pp. 1436--1445 (2017)

\bibitem{bouchitte2018_vogt}
Bouchitt{\'e}, G., Fragal{\`a}, I.: {A Duality Theory for Non-convex Problems
  in the Calculus of Variations}.
\newblock Arch Rational Mech Anal \textbf{229}(1), 361--415 (2018)

\bibitem{bredies2018_vogt}
Bredies, K., Holler, M., Storath, M., Weinmann, A.: {Total Generalized
  Variation for Manifold-Valued Data}.
\newblock SIAM J Imaging Sci \textbf{11}, 1785--1848 (2018)

\bibitem{calinescu1998_vogt}
C{\v{a}}linescu, G., Karloff, H., Rabani, Y.: {An Improved Approximation
  Algorithm for Multiway Cut}.
\newblock In: Proc STOC 1998, pp. 48--52 (1998)

\bibitem{chambolle2012_vogt}
Chambolle, A., Cremers, D., Pock, T.: {A convex approach to minimal
  partitions}.
\newblock SIAM J Imaging Sci \textbf{5}(4), 1113--1158 (2012)

\bibitem{chambolle2011_vogt}
Chambolle, A., Pock, T.: {A first-order primal-dual algorithm for convex
  problems with applications to imaging}.
\newblock {J Math Imaging Vis} \textbf{40}(1), 120--145 (2011)

\bibitem{chan2006_vogt}
Chan, T.F., Esedoglu, S., Nikolova, M.: {Algorithms for Finding Global
  Minimizers of Image Segmentation and Denoising Models}.
\newblock SIAM J Appl Math \textbf{66}, 1632--1648 (2006)

\bibitem{chan2001_vogt}
Chan, T.F., Kang, S.H., Shen, J.: {Total Variation Denoising and Enhancement of
  Color Images Based on the CB and HSV Color Models}.
\newblock J Vis Commun Image Represent \textbf{12}, 422--435 (2001)

\bibitem{chefdhotel2004_vogt}
Chefd'Hotel, C., Tschumperl{\'e}, D., Deriche, R., Faugeras, O.D.:
  {Regularizing Flows for Constrained Matrix-Valued Images}.
\newblock J Math Imaging Vis \textbf{20}, 147--162 (2004)

\bibitem{cremers2012_vogt}
Cremers, D., Strekalovskiy, E.: {Total Cyclic Variation and Generalizations}.
\newblock J Math Imaging Vis \textbf{47}, 258--277 (2012)

\bibitem{delaunoy2009_vogt}
Delaunoy, A., Fundana, K., Prados, E., Heyden, A.: {Convex multi-region
  segmentation on manifolds}.
\newblock In: Proc ICCV 2009, pp. 662--669 (2009)

\bibitem{dziuk2013_vogt}
{Dziuk}, G., {Elliott}, C.M.: {Finite element methods for surface PDEs.}
\newblock {Acta Numerica} \textbf{22}, 289--396 (2013)

\bibitem{federer1974_vogt}
{Federer}, H.: {Real flat chains, cochains and variational problems.}
\newblock {Indiana U Math J} \textbf{24}, 351--407 (1974)

\bibitem{fletcher2012_vogt}
Fletcher, P.T.: {Geodesic Regression and the Theory of Least Squares on
  Riemannian Manifolds}.
\newblock Int J Comput Vis \textbf{105}, 171--185 (2012)

\bibitem{gesch2009_vogt}
Gesch, D., Evans, G., Mauck, J., Hutchinson, J., Carswell~Jr, W.J., et~al.:
  {The national map - Elevation}.
\newblock US geological survey fact sheet \textbf{3053}(4) (2009)

\bibitem{giaquinta1998_vogt}
{Giaquinta}, M., {Modica}, G., {Sou\v{c}ek}, J.: {Cartesian currents in the
  calculus of variations I and II.}
\newblock Berlin: Springer (1998)

\bibitem{goldluecke2010_vogt}
Goldl{\"u}cke, B., Cremers, D.: {Convex Relaxation for Multilabel Problems with
  Product Label Spaces}.
\newblock In: {Proc ECCV 2010}, pp. 225--238 (2010)

\bibitem{goldluecke2013_vogt}
Goldl{\"u}cke, B., Strekalovskiy, E., Cremers, D.: {Tight Convex Relaxations
  for Vector-Valued Labeling}.
\newblock SIAM J Imaging Sci \textbf{6}, 1626--1664 (2013)

\bibitem{goldstein2013_vogt}
Goldstein, T., Esser, E., Baraniuk, R.: {Adaptive primal dual optimization for
  image processing and learning}.
\newblock In: {Proc 6th NIPS Workshop Optim Mach Learn}, pp. 1--5 (2013)

\bibitem{greig1989_vogt}
Greig, D.M., Porteous, B.T., Seheult, A.H.: {Exact Maximum A Posteriori
  Estimation for Binary Images}.
\newblock J R Stat Soc Series B Stat Methodol \textbf{51}(2), 271--279 (1989)

\bibitem{ishikawa2003_vogt}
Ishikawa, H.: {Exact Optimization for Markov Random Fields with Convex Priors}.
\newblock IEEE Trans Pattern Anal Mach Intell \textbf{25}, 1333--1336 (2003)

\bibitem{karcher1977_vogt}
{Karcher}, H.: {Riemannian center of mass and mollifier smoothing.}
\newblock {Commun Pure Appl Math} \textbf{30}, 509--541 (1977)

\bibitem{kleinberg2002_vogt}
Kleinberg, J.M., Tardos, {\'E}.: {Approximation algorithms for classification
  problems with pairwise relationships: metric labeling and Markov random
  fields}.
\newblock J ACM \textbf{49}, 616--639 (2002)

\bibitem{klodt2008_vogt}
Klodt, M., Schoenemann, T., Kolev, K., Schikora, M., Cremers, D.: {An
  Experimental Comparison of Discrete and Continuous Shape Optimization
  Methods}.
\newblock In: {Proc ECCV 2008}, pp. 332--345 (2008)

\bibitem{laude2016_vogt}
Laude, E., M{\"o}llenhoff, T., Moeller, M., Lellmann, J., Cremers, D.:
  {Sublabel-Accurate Convex Relaxation of Vectorial Multilabel Energies}.
\newblock In: {Proc ECCV 2016}, pp. 614--627 (2016)

\bibitem{laus2017_vogt}
Laus, F., Persch, J., Steidl, G.: {A Nonlocal Denoising Algorithm for
  Manifold-Valued Images Using Second Order Statistics}.
\newblock SIAM J Imaging Sci \textbf{10}, 416--448 (2017)

\bibitem{lavenant2017_vogt}
Lavenant, H.: {Harmonic mappings valued in the Wasserstein space}.
\newblock Tech. Rep. arXiv:1712.07528, arXiv (2017)

\bibitem{lee2013_vogt}
{Lee}, J.M.: {Introduction to smooth manifolds. 2nd revised ed.}, vol. 218, 2nd
  revised ed edn.
\newblock New York, NY: Springer (2013)

\bibitem{lellmann2011phd_vogt}
Lellmann, J.: {Nonsmooth convex variational approaches to image analysis}.
\newblock Ph.D. thesis, Ruprecht-Karls-Universit{\"a}t Heidelberg (2011)

\bibitem{lellmann2009_vogt}
Lellmann, J., Becker, F., Schn{\"o}rr, C.: {Convex optimization for multi-class
  image labeling with a novel family of total variation based regularizers}.
\newblock In: Proc ICCV 2009, pp. 646--653 (2009)

\bibitem{lellmann2013b_vogt}
Lellmann, J., Lellmann, B., Widmann, F., Schn\"{o}rr, C.: {Discrete and
  Continuous Models for Partitioning Problems}.
\newblock Int J Comput Vis \textbf{104}(3), 241--269 (2013)

\bibitem{lellmann2012_vogt}
Lellmann, J., Lenzen, F., Schn{\"o}rr, C.: {Optimality Bounds for a Variational
  Relaxation of the Image Partitioning Problem}.
\newblock J Math Imaging Vis \textbf{47}, 239--257 (2012)

\bibitem{lellmann2011c_vogt}
Lellmann, J., Schn{\"o}rr, C.: {Continuous Multiclass Labeling Approaches and
  Algorithms}.
\newblock SIAM J Imaging Sci \textbf{4}(4), 1049--1096 (2011)

\bibitem{lellmann2013_vogt}
Lellmann, J., Strekalovskiy, E., Koetter, S., Cremers, D.: {Total Variation
  Regularization for Functions with Values in a Manifold}.
\newblock In: Proc ICCV 2013, pp. 2944--2951 (2013)

\bibitem{loewenhauser2018_vogt}
Loewenhauser, B., Lellmann, J.: {Functional Lifting for Variational Problems
  with Higher-Order Regularization}.
\newblock In: X.C. Tai, E.~Bae, M.~Lysaker (eds.) {Imaging, Vision and Learning
  Based on Optimization and PDEs}, pp. 101--120. Springer International
  Publishing, Cham (2018)

\bibitem{massonnet1998_vogt}
Massonnet, D., Feigl, K.L.: {Radar interferometry and its application to
  changes in the Earth's surface}.
\newblock Rev Geophys \textbf{36}(4), 441--500 (1998)

\bibitem{mollenhoff2017_vogt}
M{\"o}llenhoff, T., Cremers, D.: {Sublabel-Accurate Discretization of Nonconvex
  Free-Discontinuity Problems}.
\newblock In: Proc ICCV 2017, pp. 1192--1200 (2017)

\bibitem{mollenhoff2019_vogt}
M{\"o}llenhoff, T., Cremers, D.: {Lifting Vectorial Variational Problems: A
  Natural Formulation based on Geometric Measure Theory and Discrete Exterior
  Calculus}.
\newblock In: Proc CVPR 2019 (2019)

\bibitem{mollenhoff2016_vogt}
M{\"o}llenhoff, T., Laude, E., Moeller, M., Lellmann, J., Cremers, D.:
  {Sublabel-Accurate Relaxation of Nonconvex Energies}.
\newblock In: Proc CVPR 2016 (2016)

\bibitem{pock2009_vogt}
Pock, T., Cremers, D., Bischof, H., Chambolle, A.: {An algorithm for minimizing
  the Mumford-Shah functional}.
\newblock Proc ICCV 2009 pp. 1133--1140 (2009)

\bibitem{pock2010_vogt}
Pock, T., Cremers, D., Bischof, H., Chambolle, A.: {Global Solutions of
  Variational Models with Convex Regularization}.
\newblock SIAM J Imaging Sci \textbf{3}(4), 1122--1145 (2010)

\bibitem{pock2008_vogt}
Pock, T., Schoenemann, T., Graber, G., Bischof, H., Cremers, D.: {A Convex
  Formulation of Continuous Multi-label Problems}.
\newblock In: {Proc ECCV 2008}, pp. 792--805 (2008)

\bibitem{ranftl2013_vogt}
Ranftl, R., Pock, T., Bischof, H.: {Minimizing TGV-Based Variational Models
  with Non-convex Data Terms}.
\newblock In: A.~Kuijper, K.~Bredies, T.~Pock, H.~Bischof (eds.) Proc SSVM
  2013, pp. 282--293. Springer Berlin Heidelberg, Berlin, Heidelberg (2013)

\bibitem{rocca1997_vogt}
Rocca, F., Prati, C., Ferretti, A.: {An overview of {SAR} interferometry}.
\newblock In: {Proc 3rd ERS Symp Spac Serv Env} (1997).
\newblock
  \urlprefix\url{http://earth.esa.int/workshops/ers97/program-details/speeches/rocca-et-al}

\bibitem{rosman2011_vogt}
Rosman, G., Bronstein, M.M., Bronstein, A.M., Wolf, A., Kimmel, R.:
  {Group-Valued Regularization Framework for Motion Segmentation of Dynamic
  Non-rigid Shapes}.
\newblock In: A.M. Bruckstein, B.M. ter Haar~Romeny, A.M. Bronstein, M.M.
  Bronstein (eds.) Proc SSVM 2011, pp. 725--736. Springer Berlin Heidelberg,
  Berlin, Heidelberg (2012)

\bibitem{storath2018_vogt}
Storath, M., Weinmann, A.: {Wavelet Sparse Regularization for Manifold-Valued
  Data}.
\newblock Tech. Rep. arXiv:1808.00505, arXiv (2018)

\bibitem{strecke2019_vogt}
Strecke, M., Goldluecke, B.: {Sublabel-Accurate Convex Relaxation with Total
  Generalized Variation Regularization}.
\newblock In: T.~Brox, A.~Bruhn, M.~Fritz (eds.) Proc GCPR 2018, pp. 263--277.
  Springer International Publishing (2019)

\bibitem{strekalovskiy2015_vogt}
Strekalovskiy, E.: {Convex Relaxation of Variational Models with Applications
  in Image Analysis}.
\newblock Ph.D. thesis, Technische Universit{\"a}t M{\"u}nchen (2015)

\bibitem{strekalovskiy2011b_vogt}
Strekalovskiy, E., Cremers, D.: {Total variation for cyclic structures: Convex
  relaxation and efficient minimization}.
\newblock In: Proc CVPR 2011, pp. 1905--1911 (2011)

\bibitem{strekalovskiy2011_vogt}
Strekalovskiy, E., Goldl{\"u}cke, B., Cremers, D.: {Tight convex relaxations
  for vector-valued labeling problems}.
\newblock Proc ICCV 2011 pp. 2328--2335 (2011)

\bibitem{strekalovskiy2012_vogt}
Strekalovskiy, E., Nieuwenhuis, C., Cremers, D.: {Nonmetric Priors for
  Continuous Multilabel Optimization}.
\newblock In: A.~Fitzgibbon, S.~Lazebnik, P.~Perona, Y.~Sato, C.~Schmid (eds.)
  {Proc ECCV 2012}, pp. 208--221. Springer Berlin Heidelberg (2012)

\bibitem{vogt2018_vogt}
Vogt, T., Lellmann, J.: {Measure-Valued Variational Models with Applications to
  Diffusion-Weighted Imaging}.
\newblock J Math Imaging Vis \textbf{60}, 1482--1502 (2018)

\bibitem{vogt2019_vogt}
Vogt, T., Lellmann, J.: {Functional Liftings of Vectorial Variational Problems
  with Laplacian Regularization}.
\newblock In: M.~Burger, J.~Lellmann, J.~Modersitzki (eds.) Proc SSVM 2019, pp.
  559--571 (2019)

\bibitem{weinmann2014_vogt}
Weinmann, A., Demaret, L., Storath, M.: {Total variation regularization for
  manifold-valued data}.
\newblock SIAM J Imaging Sci \textbf{7}, 2226--2257 (2014)

\bibitem{weinmann2015_vogt}
Weinmann, A., Demaret, L., Storath, M.: {Mumford-Shah and Potts Regularization
  for Manifold-Valued Data}.
\newblock J Math Imaging Vis \textbf{55}, 428--445 (2015)

\bibitem{windheuser2016_vogt}
Windheuser, T., Cremers, D.: {A Convex Solution to Spatially-Regularized
  Correspondence Problems}.
\newblock In: B.~Leibe, J.~Matas, N.~Sebe, M.~Welling (eds.) {Proc ECCV 2016},
  pp. 853--868. Springer International Publishing (2016)

\bibitem{zach2008_vogt}
Zach, C., Gallup, D., Frahm, J.M., Niethammer, M.: {Fast Global Labeling for
  Real-Time Stereo Using Multiple Plane Sweeps}.
\newblock In: Proc VMV 2008, pp. 243--252 (2008)

\bibitem{zach2012_vogt}
Zach, C., Kohli, P.: {A Convex Discrete-Continuous Approach for Markov Random
  Fields}.
\newblock In: {Proc ECCV 2012}, pp. 386--399 (2012)

\end{thebibliography}
